\documentclass{article}
\pdfoutput=1
\usepackage{makeidx}        
\usepackage{graphicx}      
\usepackage[bottom]{footmisc}

\usepackage{amssymb,latexsym,amsmath, amsthm}

\newcommand{\dwedge}[2]%
{\overset{#1}{\underset{#2}{\mbox{$\bigwedge\hspace{-2ex}\bigwedge$}}}}
\newcommand{\dvee}[2]%
{\overset{#1}{\underset{#2}{\mbox{$\bigvee\hspace{-2ex}\bigvee$}}}}

\newcommand{\inMath}[1]{\relax\ifmmode{#1}%
\else{\mbox{$#1$}}\fi}

\newcommand{\vdot}%
{\unitlength0.4mm
\begin{picture}(8,10)
\put(4,0){.}
\put(0,10){.}
\put(8,10){.}
\put(2,5){.}
\put(6,5){.}
\end{picture}
}

\newcommand{\ignore}[1]{}

\newcommand{\romanref}[1]{%
\if \ref{#1}\empty {\setcounter{romanrefcounter}0} \else 
{\setcounter{romanrefcounter}{\ref{#1}}}\fi%
{\it \roman{romanrefcounter}}}

\newcommand{\Bc}{{\bf c}}

\newcommand{\BC}{{\bf C}}

\newcommand{\BD}{{\bf D}}

\newcommand{\Bd}{{\bf d}}

\newcommand{\BI}{{\bf I}}

\newcommand{\BK}{{\bf K}}

\newcommand{\BL}{{{\bf L}}}

\newcommand{\BQ}{{\bf Q}}

\newcommand{\BR}{{\bf R}}

\newcommand{\BS}{{\bf S}}

\newcommand{\BU}{{\bf U}}

\newlength{\circlen}
\newlength{\symblen}

\ifx\text\undefined
  \newcommand{\text}[1]{\relax
    \ifmmode\mathchoice
      {\hbox{\the\textfont0\relax#1}}%
      {\hbox{\the\textfont0\relax#1}}%
      {\hbox{\the\scriptfont0\relax#1}}%
      {\hbox{\the\scriptscriptfont0\relax#1}}%
    \else{\relax#1}\fi}
  \fi

\newcommand{\defsub}[1]{}

\newcommand{\nin}{{\not\in}}

\def\twoheaddownarrow{\rlap{$\downarrow$}\raise-.5ex\hbox{$\downarrow$}}
\def\twoheaduparrow{\rlap{$\uparrow$}\raise.5ex\hbox{$\uparrow$}}

\def\texturespicture #1 by #2 (#3){
\vbox to #2 {\hrule width #1 height 0pt depth 0pt
}}

\def\scaledpicture #1 by #2 (#3 scaled #4){{\dimen0=#1 \dimen1=#2
\divide\dimen0 by 1000 \multiply \dimen0 by #4
\divide\dimen1 by 1000 \multiply \dimen1 by #4
\texturespicture\dimen0 by \dimen1 (#3 scaled #4)}}

 {\begin{tabular}[t]{llr} && \fbox{#1}\\}%
 {\end{tabular}}

{\end{list}}

{\end{tabular}%
\end{trivlist}}

\newcommand\strikethrough[1]{{\setbox0=\hbox{$#1$}
\hrule height.75ex depth-.65ex width\wd0 \kern-\wd0\box0}}

\mathchardef\gt="313E \mathchardef\lt="313C

\def\undern#1{\vtop{\ialign{##\crcr
 $\hfil\displaystyle{#1}\hfil$\crcr
 \noalign{\kern-.1pt\nointerlineskip}%
 ~\,\raisebox{-.5ex}{$n \to \infty$} \crcr}}}

\def\undern#1{\vtop{\ialign{##\crcr
 $\hfil\displaystyle{#1}\hfil$\crcr
 \noalign{\kern-.1pt\nointerlineskip}%
 ~\,\raisebox{-.5ex}{$i$} \crcr}}}

\newtheorem{theorem}{Theorem}[section]

\renewcommand{\labelenumi}{(\roman{enumi})}

\newfont{\bi}{cmmib10}
\newtheorem{coro}[theorem]{Corollary}
\newtheorem{lem}[theorem]{Lemma}
\newtheorem{propo}[theorem]{Proposition}

\newtheorem{rem}[theorem]{Remark}

\newtheorem{defi}[theorem]{Definition}

\def\<{\langle}
\def\>{\rangle}

\def\Lrarr{{\Leftrightarrow}}

\def\nvd{{\not\vdash}}

\def\a1an{{a_1,\dots,a_n}}

\def\al{\alpha}
\def\al{\alpha}

\def\conj{~\&~}

\def\b1{\Box^{-1}}
\def\B1Bk{{B_1,\dots,B_k}}

\def\baF{\overline{F}}
\def\baR{\overline{R}}

\def\bbK4{{\bf K4}}
\def\bbD4{{\bf D4}}
\def\bbS4{{\bf S4}}

\def\bD{{\bf D}}
\def\bD4{{\bf D4}}

\def\be{\beta}

\def\bK{{\bf K}}
\def\bL{{\tx{\bf {\em L}}}}
\def\bLa{\mathbf{\Lambda}}

\def\bQS4{{\bf QS4}}

\def\bS5{{\bf S5}}
\def\BS5{{\bf S5}}

\def\bvee0{\bigvee}
\def\bvee{\boldsymbol{\vee}}

\def\C12{{\cal C}_1\t{\cal C}_2}

\def\cC{{\cal C}}

\def\cd{\cdot}
\def\cD{{\cal D}}

\def\cL{{\cal L}}

\def\cP{{\cal P}}

\def\d1{\Diamond^{-1}}
\def\da1aian{{a_1,\dots,a_i,\dots,a_n}}
\def\da1ban{{a_1,\dots,b,\dots,a_n}}
\def\da1n{{a_1,\dots,a_n}}

\def\Del{{\Delta}}

\def\Di{\Diamond}

\def\dri{\twoheadrightarrow}

\def\dS5{{\text{{\bf {\em S5}}}}}

\def\dx1n{\tx{\hbox{$x_1,\dots,x_n$}}}
\def\dx1yxn{{x_1,\dots,y,\dots,x_n}}

\def\emp{\varnothing}

\def\eqv{\equiv}

\def\ex{\exists}

\def\fo{\forall}

\def\gam{\gamma}

\def\I{\Phi}

\def\lora{\longrightarrow}
\def\Lra{\Leftrightarrow}

\def\lrarr{\leftrightarrow}

\def\Lri{\Longleftrightarrow}

\def\miff{\;\mbox{iff}\;}

\def\mo{\vDash}

\def\mt{\mapsto}

\def\nin{\not\in}
\def\nmo{{\,\not\mo\,}}

\def\oB{\overline{\square}}

\def\om{\omega}

\def\p1pk{{p_1,\dots,p_k}}

\def\r1rn{{R_1,\dots,R_n}}
\def\r1rN{{R_1,\dots,R_N}}
\def\Ra{\Rightarrow}

\def\ri{\rightarrow}
\def\rr1n{{r_1,\dots,r_n}}

\def\s1sn{{S_1,\dots,S_n}}
\def\sbq{\subseteq}

\def\sp{\supset}

\def\spq{\supseteq}

\def\sqr{\square}
\def\ssb{\subset}

\def\t{\times}

\def\tF{\tilde{F}}

\def\th{\theta}

\def\ti{\times}

\def\tW{\tilde{W}}
\def\tx{\text}

\def\uR{{\underline{R}}}

\def\vd{\vdash}

\def\vp{\varphi}

\def\we{\wedge}

\def\x1xn{{x_1,\dots,x_n}}
\def\xx1n{{x,x_1,\dots,x_n}}

\def\l1{\L}

\def\Ri0{{R_{i_0}}}

\def\Rj0{{R_{j_0}}}

\def\Rn{{\bf R}^n}

\def\vapi0{{\varphi_{i_0}}}

\def\vapj0{{\varphi_{j_0}}}

\def\cL{{\cal L}}

\def\emp{{\varnothing}}

\def\lrarr{{\leftrightarrow}}

\def\r1rN{{R_1,\dots,R_N}}

\newcommand{\X}{\mathfrak X}
\newcommand{\YY}{\mathfrak Y}
\newcommand{\DA}{{[\ne]}}
\newcommand{\DE}{{\langle\ne\rangle}}
\newcommand{\UA}{\left[\forall\right]}

\renewcommand{\I}{{\bf I}}

\newcommand{\romb}{\Diamond}

\newcommand{\pmor}{\twoheadrightarrow}
\newcommand{\dpmor}{\pmor^d{}}

\newcommand{\norm}[1]{\left\Vert#1\right\Vert}
\newcommand{\abs}[1]{\left\vert#1\right\vert}
\newcommand{\set}[1]{\left\{#1\right\}}

\newcommand{\setdef}[2][x]{\left \{#1\,\left|\,#2 \right. \right \}}
\newcommand{\Real}{\mathbf R}

\newcommand{\Y}{\mathcal Y}
\newcommand{\Z}{\mathcal Z}

\newcommand{\cpmor}{\twoheadrightarrow^{dd}}

\newcommand{\impl}{\supset}

\renewcommand{\bLa}{\Lambda}
\renewcommand{\sp}{\to}

\newcommand{\lgc}[1]{\mathbf{#1}}

\newtheorem*{prop*}{Proposition}{\bfseries}{\itshape}

\title{Derivational modal logics with the difference modality}
\author{Andrey Kudinov\\ \small
kudinov\ [at here]\ iitp\ [dot]\ ru \\
\small
Institute for Information Transmission Problems, Russian Academy of
Sciences\\ \small
National Research University Higher School of Economics, Moscow, Russia\\ \small
Moscow Institute of Physics and Technology\\ \and
Valentin Shehtman \\ \small shehtman\ [at here]\ netscape\ [dot]\ net \\
\small Institute for Information Transmission Problems, Russian Academy of Sciences\\ \small
Moscow State University\\ \small
National Research University Higher School of Economics, Moscow, Russia}

\begin{document}

\maketitle

\abstract{In this chapter we study modal logics of topological spaces in the combined language with the
derivational modality and the difference modality. We give 
axiomatizations and prove completeness for the following classes: all spaces, $T_1$-spaces, dense-in-themselves spaces,
a zero-dimensional dense-in-itself separable metric space, 
$\Real^n~
(n\ge 2)$. We also discuss the correlation between languages with different combinations of the topological,
the derivational, the universal and the difference modality in terms of definability.}


\section{Introduction}
\label{sec:intro} Topological modal logic was initiated by the works of A. Tarski and J.C.C. McKinsey in the 1940s. 
They were first to consider both topological interpretations of the diamond modality: one as closure, and another as derivative. 
	
Their studies of closure modal 
logics were rather detailed and profound. 
In particular, in the fundamental paper \cite{MT} they have shown that the logic of any metric separable dense-in-itself space is $\bbS4$. 
This remarkable result also demonstrates a relative weakness of the closure operator to distinguish between interesting topological properties. 
       
The derivational interpretation gives more expressive power. For example, the real line can be distinguished from the real plane (the observation made by K. Kuratowski as early as in 1920s, cf. \cite{K22});  the real line can be distinguished from the rational line \cite{Sh90}; $T_0$ and $T_D$ separation axioms become expressible \cite{BEG}, \cite{E1}. However, in \cite{MT} McKinsey and Tarski only gave basic definitions for derivational modal logics and put several problems that were solved much later.
       
The derivational semantics also has its limitations (for example, it is still impossible to distinguish $\BR^2$ from $\BR^3$). Further increase of expressive power can be provided by the well-known methods of adding universal or difference modalities \cite{GorPas}, \cite{GarGor}.  In the context of topological semantics this approach also has proved fruitful --- for example, connectedness is expressible in modal logic with the closure and the universal modality \cite{Sh99}, and the $T_1$ separation axiom in modal logic with the closure and the difference modality \cite{kudinov_2005}.  
       
Until the early 1990s, when the connections between topological modal logic and Computer Science were established, the interest in that subject was moderate. Leo Esakia was one of the enthusiasts of modal logical approach to topology, and he was probably the
first to appreciate the role of the derivational modality, in particular, in modal logics of provability \cite{E}. Another strong motivation for further studies of derivational modal logics (`d-logics') were the axiomatization problems left open in \cite{MT}.\footnote{The early works of the second author in this field were greatly influenced by Leo Esakia.} In recent years d-logics have been studied rather intensively, a brief summary of results can be found in section 3 below.

In this chapter the first thorough investigation is provided for logics in the most expressive language in this context\footnote{Some other kinds of topomodal logics arise when we deal with topological spaces with additional structures, e.g. spaces with two topologies, spaces with a homeomorphism etc. (cf. \cite{H}).}, 
namely the derivational modal logics with the difference modality (`dd-logics'). It unifies earlier studies by the first author in closure modal logics with the difference modality (`cd-logics') and by the second author in d-logics.

The diagram in section 12 compares the expressive power of different kinds of topomodal logics. Our conjecture 
is that dd-logics are strictly more expressive than the others, but it is still an open question if the dd-language is stronger than the cd-language. Speaking informally, it is more convenient --- for example, the Kuratowski's axiom for $\BR^2$ (Definition 9.1) is expressible in cd-logic as well, but in a more complicated form \cite{kudinov_2006}.

We show that still in many cases 
properties of dd-logics are similar to those of d-logics: finite axiomatizability, decidability and the
finite model property (fmp).
Besides specific results characterizing logics of some particular spaces, our goal was to propose some general methods. In fact, nowadays in topomodal logic there are many technical proofs, but few general methods. In this chapter we propose only two simplifying novelties --- dd-morphisms (section 6) and the Glueing lemma \ref{Glue}, but we hope that much more can be done in this direction, cf. the recent paper \cite{Hodkinson}.

In more detail, the plan of the chapter is as follows. Preliminary sections 2--4 include standard
definitions and basic facts about modal logics and their semantics. Some general completeness results
for dd-logics can be found in sections 5, 7. In section 5 we show that every extension of the minimal
logic $\lgc{K4^\circ D^+}$ by variable-free axioms is topologically complete. In section 8 we prove
the same for extensions of $\lgc{DT_1}$ (the logic of dense-in-themselves $T_1$-spaces); the proof is
based on a construction of d-morphisms from the recent paper \cite{BLB}.

In section 6  we consider validity-preserving maps from topological to Kripke frames (d-morphisms and
dd-morphisms) and prove a modified version of McKinsey--Tarski's lemma on dissectable spaces. In
section 7 we prove that $\lgc{DT_1}$ is complete w.r.t. an arbitrary zero-dimensional dense-in-itself
separable metric space by the method from \cite{Sh90}, \cite{Sh00}.

Sections 8--10 study the axiom of connectedness $AC$ and Kuratowski's axiom $Ku$ related to local
1-componency. In particular we prove that the logic $\lgc{DT_1CK}$ with both these axioms has the
fmp. This is a refinement of an earlier result \cite{Sh90}, \cite{Sh00} on the fmp of the d-logic
$\lgc{D4}+Ku$ (the new proof uses a simpler construction).

Section 11 contains our central result: $\lgc{DT_1CK}$ is the dd-logic of $\BR^n$ for $n>1$. The
proof uses an inductive construction of dd-morphisms onto finite frames of the corresponding logic,
and it combines methods from \cite{Sh90}, \cite{Sh00}, \cite{kudinov_2006}, with an essential
improvement motivated by \cite{LB2} and based on the Glueing lemma. 

The final section discusses some further directions and open questions. The Appendix contains technical details of some proofs.

\section{Basic notions}
\label{sec:bas}
\renewcommand{\labelenumi}{(\arabic{enumi})}

The material of this section is quite standard, and most of it can be found in \cite{CZ97}.
We
consider {\em $n$-modal (propositional) formulas} constructed from a countable set of propositional
variables $PV$ and the connectives $\bot$, $\ri$, $\sqr_1,\ldots,\sqr_n$. The derived connectives are
$\wedge,\;\vee,\;\neg,\;\top,\;\lrarr,\;\Di_1,\ldots,\Di_n$. A formula without occurrences of
propositional variables is called  {\em closed}.

A {\em  (normal) $n$-modal logic } is a set of modal formulas containing the classical tautologies,
the axioms $\sqr_i(p\ri q)\ri (\sqr_i p \ri\sqr_i q)$ and closed under the standard inference rules:
Modus Ponens ($A,~A\ri B/B$), Necessitation ($A/\sqr_i A$), and Substitution $(A(p_j)/A(B))$.

To be more specific, we use the terms `($\sqr_1,\ldots,\sqr_n$)-modal formula' and
`($\sqr_1,\ldots,\sqr_n$)-modal logic'.

$\BK_n$ denotes the minimal $n$-modal logic (and $\BK=\BK_1$). An $n$-modal logic containing a
certain $n$-modal logic $\bLa$ is called an \emph{extension} of $\bLa$, or a \emph{$\bLa$-logic}. The
minimal $\bLa$-logic containing a set of $n$-modal formulas $\Gamma$ is denoted by $\bLa+\Gamma$. In
particular,
$$\bbK4:=\BK+\sqr p\ri\sqr\sqr p,~
\bbS4:=\bbK4+\sqr p\ri p, ~\bbD4:=\bbK4+\Di\top,
$$
$$
\bbK4^\circ := {\bf wK4}:= \bK + p\we \sqr p\ri \sqr \sqr p.
$$

The {\em fusion} $L_1*L_2$ of modal logics $L_1,~L_2$ with distinct modalities is the smallest modal
logic in the joined language containing $L_1\cup L_2$.

A {\em (normal)  $n$-modal algebra} is a Boolean algebra with extra $n$ unary operations preserving
${\bf 1}$ (the unit) and distributing over $\cap$;  they are often denoted by $\sqr_1,\ldots,\sqr_n$,
in the same way as the modal connectives. A {\em valuation} in a    modal algebra $\mathfrak{A}$ is a
set-theoretic map $\th: PV \lora \mathfrak{A}$. It extends to all $n$-modal formulas by induction:
$$\th(\bot)=\emp,~\th(A\ri B)=-\th(A)\cup\th(B),~
\th(\sqr_i A)=\sqr_i\th(A).$$ A formula $A$ is {\em true in} $\mathfrak{A}$ (in symbols:
$\mathfrak{A}\mo A$) if $\th(A)=\textbf{1}$  for any valuation $\th$. The set $\BL(\mathfrak{A})$ of
all $n$-modal formulas true in an $n$-modal algebra $\mathfrak{A}$ is an $n$-modal logic called the
{\em  logic of $\mathfrak{A}$}.

An $n$-modal  {\em Kripke frame } is a tuple $ F=(W,R_1,\ldots,R_n)$, where $W$ is a nonempty set
(of worlds), $R_i$ are binary relations on $W$. We often write $x\in F$ instead of $x\in W$. In this
chapter (except for Section 2) all 1-modal frames are assumed to be transitive. The associated $n$-modal algebra $MA(F)$ is $2^W$ (the Boolean
algebra of all subsets of $W$) with the operations $\sqr_1,\ldots,\sqr_n$ such that $\sqr_i V= \{x\mid R_i(x)\sbq
V\}$ for any $V\sbq W$.

A {\em valuation} in $F$ is the same as in $MA(F)$, i.e., this is   a map from $PV$ to $\cP(W)$ (the
power set of $W$). A {\em (Kripke) model} over $F$ is a pair $M=(F,\th)$, where $\th$ is a valuation
in $F$. The notation $M,x\mo A$ means $x\in\th(A)$, which is also read as `$A$ is \emph{true in $M$
at} $x$'. A (modal) formula $A$ is {\em true in} $M$ (in symbols:  $M\mo A$) if $A$ is true in $M$ at
all worlds. A formula $A$ is called {\em valid in} a Kripke frame $F$ (in symbols: $F\mo A$) if $A$
is true in all Kripke models over $F$; this is obviously equivalent to $MA(F)\mo A$.


The {\em modal logic} $\BL(F)$ of a Kripke frame $F$ is the set of all modal formulas valid in $F$,
i.e., $\BL(MA(F))$. For a class of $n$-modal frames $\cC$, the {\em modal logic} of $\cC$ (or \emph{the modal logic determined} by $\cC$) is $\BL(\cC):=\bigcap\{\BL(F)\mid F\in\cC\}$. Logics determined by classes of
Kripke frames are called {\em Kripke complete}. An $n$-modal frame validating an $n$-modal logic
$\bLa$ is called a {\em $\bLa$-frame}. A modal logic has the {\em finite model property (fmp)} if it
is determined by some class of finite frames.

It is well known that $(W,R)\mo\bbK4$ iff $R$ is transitive;
$(W,R)\mo\bbS4$ iff $R$ is reflexive transitive (a {\em quasi-order}).

A {\em  cluster} in a transitive frame $(W,R)$ is an equivalence class under the relation $\sim_R:=(R
\cap R^{-1}) \cup I_W$, where $I_W$ is the equality relation on $W$. A {\em degenerate cluster} is an
irreflexive singleton. A cluster that is a reflexive singleton, is called  {\em trivial}, or {\em simple}. A {\em chain}  is a frame $(W,R)$ with $R$
transitive, antisymmetric and linear, i.e., it satisfies $\fo x\fo y~
(xRy\vee yRx\vee x=y)$. A point $x\in W$ is {\em strictly (R-)minimal} if $R^{-1}(x)=\emp$.

A {\em subframe} of a frame $F=(W, R_1,\ldots,R_n)$ obtained by restriction to $V\sbq W$, is
$F|V:=(V,R_1|V,\ldots,R_n|V)$. Then for any Kripke model $M=(F,\th)$ we have a {\em submodel}
$M|V:=(F|V,\th|V)$, where $(\th|V)(q):=\th(q)\cap V$ for each $q\in PV$. If $R_i(V)\sbq V$ for any
$i$, the subframe $F|V$ and the submodel $M|V$ are called {\em generated}.

The {\em union} of subframes $F_j=F|W_j,~j\in J$ is the subframe
$\bigcup\limits_{j\in J}F_j:=F|\bigcup\limits_{j\in J}W_j$.

A  {\em generated subframe (cone) with the root $x$} is
$F^x:=F|R^*(x)$, where $R^*$ is the reflexive transitive closure
of $R_1\cup\ldots\cup R_n$; so for a transitive frame $(W,R)$,~$R^*=R\cup I_W$ is the reflexive
closure of $R$ (which is also denoted by $\baR$). A frame $F$ is called {\em rooted} with the root
$u$ if $F=F^u$. Similarly we define a cone $M^x$ of a Kripke model $M$.

Every finite rooted  transitive frame $F=(W,R)$ can be presented as the union $(F|C)\cup
F^{x_1}\cup\ldots\cup F^{x_m}$ ($m\geq 0$), where $C$ is the root cluster, $x_i$ are its successors
(i.e., $x_i\nin C,~\overline{R}^{-1}(x_i)=\sim_R(x_i)\cup C$). If $C$ is non-degenerate, the frame
$F|C$ is $(C,C^2)$, which we usually denote just by $C$. If $C=\{a\}$ is degenerate, $ F|C$ is
$(\{a\},\emp)$, which we denote by $\breve{a}$.

Let us fix the propositional language (and the number $n$) until the end of this section.
\begin{lem}\label{gl}  {\em  (Generation Lemma)}
\begin{enumerate}
\item
$\BL(F)=\bigcap\{\BL(F^x) \mid x\in F\}$.
\item
If $F$ is a generated subframe of $G$, then $\BL(G)\sbq \BL(F)$.
\item
If $M$ is a generated submodel of $N$, then for any formula $A$ for any $x$ in $M$
\[
N,x\mo A\miff M,x\mo A.
\]
\end{enumerate}
\end{lem}

\begin {lem}\label{root}
For any Kripke complete modal logic $\bLa$,
\[
\bLa=\BL(\mbox{all }\bLa\mbox{-frames})=\BL(\mbox{all rooted }\bLa\mbox{-frames}).
\]
\end{lem}

A {\em  p-morphism } from a frame $(W,R_1,\ldots,R_n)$ onto a frame $
(W^\prime,R_1^\prime,\ldots,R_n^\prime)$  is a surjective map $f:W\lora W^\prime $ satisfying the
following conditions (for any $i$):
\begin{enumerate}
\renewcommand{\labelenumi}{(\arabic{enumi})}
\item
$\fo x\fo y~ (xR_i y\Ra f(x) R_i'f(y))$ (monotonicity);
\item
$\fo x\fo z~ (f(x)R_i' z\Ra\ex y
 (f(y) = z \conj xR_iy))$ (the lift property).
\end{enumerate}
If $xR_iy$ and $f(x)R_i'f(y)$, we say that $xR_iy$ {\em lifts} $f(x)R_i'f(y)$.

Note that (1)$\conj$(2) is equivalent to
$$\fo x~f(R_i(x))=R_i'(f(x)).$$
$f:\; F\dri F'$ denotes that $f$ is
a p-morphism from $F$ onto
$F^\prime $.

Every set-theoretic map $f:W\lora W'$ gives rise to the dual morphism of Boolean algebras
$2^f:2^{W'}\lora 2^W$ sending every subset $V\sbq W'$ to its inverse image $f^{-1}(V)\sbq W$.

\begin {lem}\label{p1} {\em  (P-morphism Lemma)}~
\renewcommand{\labelenumi}{(\arabic{enumi})}
\begin{enumerate}
\item
$f:\; F\dri F'$ iff $2^f$ is an embedding of $MA(F')$ in $MA(F)$.
\item
$f:\; F\dri F'$ implies $\bL(F)\sbq\bL(F').$
\item
If $f:\; F\dri F'$, then $F\mo A\Lra F'\mo A$ for any closed formula $A$.
\end{enumerate}
\end{lem}


In proofs of the fmp in this chapter we will use the well-known filtration method \cite{CZ97}. Let us recall the construction we need.

Let   $\Psi$  be a set of modal formulas closed under subformulas.
For a
Kripke model $M =
(F,\vp )$ over a frame $F = (W,R_1,\ldots,R_n)$, there is the equivalence relation on $W$
$$ x
\eqv_\Psi y \Lri  \fo A\in\Psi  (M,x  \mo   A \Lra M,y \mo  A).$$ Put
$W^{\prime}:=W/\eqv_\Psi;~~x^\sim\,:=\; \eqv_\Psi(x)$ (the equivalence class of $x$),\\
$\vp^\prime (q)
:=\{x^\sim\mid x\in \vp (q)\}$ for $q    \in PV\cap\Psi$ (and let $\vp'(q)$ be arbitrary for $q
\in PV-\Psi$).

\begin{lem}\label{L26}{\em(Filtration Lemma)}
Under the above assumptions,
consider  the relations $\uR_i, R'_i $ on $W^\prime$ such that
$$a\uR_ib  \hbox{ iff } \ex x\in a~ \ex y\in b~ xR_iy ,$$
\[
R'_i=
\begin{cases}
\mbox{the transitive closure of }\uR_i & \mbox{ if }R_i\mbox{ is transitive,}\nonumber\\
\uR_i & \mbox{otherwise.}\nonumber
 \end{cases}
\]
Put $M^\prime := (W^\prime, R'_1,\ldots,R'_n, \vp^\prime)$.
Then
for any  $x\in W, ~A\in\Psi $  :
$$M,x \mo  A\mbox{  iff  }M',x^\sim \mo  A.$$
\end{lem}


\begin{defi}
An {\em  $m$-formula}  is a modal formula in propositional variables $\{ p_1,\dots,p_m\}$. For a
modal logic $\bLa$ we define the {\em $m$-weak} (or {\em $m$-restricted) canonical frame
$F_{\bLa\lceil m}  := (W, R_1,\ldots, R_n)$ and canonical model} $M_{\bLa\lceil m}:= (F_{\bLa\lceil
m}, \vp)$, where $W$ is the set of all maximal $\bLa$-consistent sets of $m$-formulas,
$
 x R_i y \mbox{  iff  for any }m\mbox{-formula } A$\\
$(\square_i A\in x \Ra A\in y), $
\[
\vp(p_i): =
\begin{cases}
\{ x\mid p_i\in x\}& \mbox{ if }i\leq m,\nonumber\\
\emp& \mbox{ if }i> m.
\nonumber
 \end{cases}
\]
$\bLa$ is called {\em weakly canonical} if $ F_{\bLa\lceil m}\mo\bLa$ for any finite $m$.
\end{defi}

\begin{propo}\label{cm}
For any $m$-formula $A$ and a modal logic $\bLa$
\begin{enumerate}
\item
$M_{\bLa\lceil m},x \mo  A \mbox{ iff }A\in x$;
\item
$M_{\bLa\lceil m}\mo  A  \mbox{ iff }A\in\bLa$;
\item
if $\bLa$ is weakly canonical, then it is Kripke complete.
\end{enumerate}
\end{propo}

\begin{coro}\label{cm1}
If for any $m$-formula $A$, $M_{\bLa\lceil m},x\mo A\Lrarr M_{\bLa\lceil m},y\mo A$, then $x=y$.
\end{coro}

\begin{defi}
A cluster $C$ in a transitive frame $(W,R)$ is called
\emph{maximal} if  \mbox{$\baR(C)
= C$.}
\end{defi}

\begin{lem}\label{L82}\label{L84}
Let $ F_{\bLa\lceil m} =(W,R_1,\ldots,R_n)$ and suppose $\bLa\vd\sqr_1p\ri\sqr_1\sqr_1p $ (i.e.,
$R_1$ is transitive). Then every generated subframe of $(W,R_1)$ contains a maximal cluster.
\end{lem}

The proof is based on the fact that the general Kripke frame corresponding to a canonical model is
descriptive; cf. \cite{CZ97}, \cite{F85} for further details\footnote{For the 1-modal case this lemma
has been known as folklore since the 1970s; the second author learned it from Leo Esakia in 1975.}.

\section{Derivational modal logics}
\renewcommand{\labelenumi}{(\arabic{enumi})}
We denote topological spaces by $\mathfrak{X}, \mathfrak{Y},\ldots$ and the corresponding sets by
$X,Y,\dots$.\footnote{Sometimes we neglect this difference.} The interior operation in a space $\X$
is denoted by $\BI_X$ and the closure operation by $\BC_X$, but we often omit the subscript $X$. A
set $S$ is a {\em neighbourhood} of a point $x$ if $x\in \BI S$; then $S-\{x\}$ is called a {\em
punctured neighbourhood} of $x$.

\begin{defi}
Let $\X$ be a topological space, $V\sbq X$. A point $x\in X$  is said to be {\em limit} for $V$  if
$x\in\BC(V-\{x\})$; a non-limit point of $V$ is called {\em  isolated}.

The {\em  derived  set of $V$} (denoted by $\Bd V$,
or by $\Bd_X V$) is the set of all limit points of $V$. The unary operation $V \mt \Bd V$ on
$\cP(X)$ is called {\em  the derivation} (in $\X$).

A set without isolated points is called
{\em dense-in-itself}.
\end{defi}

\begin{lem}\label{dY}\cite{K66}
For a subspace $\Y\sbq \X$ and $V\sbq X$
$\Bd_Y(V\cap Y)=\Bd_X(V\cap Y)\cap Y$; if $Y$ is open, then 
$\Bd_Y(V\cap Y)=\Bd_XV\cap Y$. 
\end{lem}

\begin{defi}
The \emph{derivational algebra of a topological space $\X$} is $DA(X):=(2^X, {\bf \tilde{d}} )$, where $2^X$ is the Boolean algebra  of all
subsets of $X$, ${\bf \tilde{d}} V := -\Bd (-V)$\footnote{There is no common notation for this operation; some authors use
$\mathbf{\tau}$.}.
The \emph{closure algebra of a space $\X$} is $CA(\X):=(2^X,\BI )$.
\end{defi}

\begin{rem}\rm
In \cite{MT} the derivational algebra of $\X$ is defined as $(2^X,\Bd )$, and the closure algebra as
$(2^X,\BC )$, but here we adopt equivalent dual definitions.
\end{rem}

It is well known that $CA(\X), ~DA(\X)$ are modal algebras, $CA(\X)\mo \bbS4$ and
$DA(\X) \mo \bbK4^\circ$ (the latter is due to Esakia).

Every Kripke $\bbS4$-frame $F=(W,R)$ is associated with a topological space $N(F)$ on $W$, with the {\em Alexandrov} (or {\em right}) {\em topology} $\{V\sbq W\mid R(V)\sbq V\}$. In $N(F)$
we have $\BC V=R^{-1}(V), ~\BI V=\{x\mid R(x)\sbq V\}$; thus $MA(F)=CA(N(F))$.

\begin{defi}
A modal formula $A$ is called {\em  d-valid in a topological space} $\X$ (in symbols, $\X\mo^d A$) if it is true  in the algebra
$DA(\X)$.  The logic $\BL(DA(\X))$ is called the {\em  derivational modal logic} (or the {\em  d-logic}) of $\X$ and denoted by $\BL\Bd (\X)$.

A formula $A$ is called {\em  c-valid} in $\X$ (in symbols, $\X\mo^c A$) if it is true in $CA(\X)$.  
$\BL\Bc (\X):=\BL(CA(\X))$ is called the {\em  closure modal logic,} or the
\emph{c-logic} of  $\X$. 
\end{defi}

\begin{defi}
For a class of topological spaces $\cC$  we also define the
{\em  d-logic}
$\BL\Bd(\cC):=\bigcap\{\BL\Bd(\X)\mid \X\in\cC\}$ and the
{\em c-logic}
$\BL\Bc(\cC):=\bigcap\{\BL\Bc(\X)\mid \X\in\cC\}$. Logics of this form
are called {\em d-complete} (respectively, {\em c-complete} ).
\end{defi}

\begin{defi}
A {\em valuation} in a topological space $\X$ is a map
$\vp : PV \lora \cP(\X)$. Then $(\X,\vp)$ is called a {\em topological model} over $\X$.
\end{defi}

So valuations in $\X$, $CA(\X)$, and $DA(\X)$ are the same. Every valuation $\vp$ can be prolonged to
all formulas in two ways, according either to $CA(\X)$ or $DA(\X)$. The corresponding maps are
denoted respectively by $\vp_c$ or $\vp_d$. Thus
\begin{align*}
\vp_d (\sqr A) &= {\bf \tilde{d}}\vp_d (A),~&
\vp_d (\Di A) &= \Bd\vp_d (A),\\
\vp_c (\sqr A) &= \BI\vp_c (A),~&
\vp_c (\Di A) &= \BC\vp_c (A). 
\end{align*}

A formula $A$ is called {\em d-true} (respectively,  {\em c-true}) in $(\X,\vp)$ if $\vp_d ( A) =X$ (respectively, $\vp_c ( A) =X$ ).
So $A$ is d-valid in $\X$ iff $A$ is d-true in every topological model over $\X$, similarly for c-validity.

\begin{defi}
A modal formula $A$ is called \emph{d-true at a point} $x$ in a topological model $(\X,\vp)$  if  $x\in\vp_d (A)$.
\end{defi}

Instead of  $x\in\vp_d(A)$, we write $x \mo^d  A$ if the model is clear from the context.
Similarly we define the c-truth at a point and use the corresponding notation.

From the definitions we obtain
\begin{lem}
For a topological model over a space $\X$
\begin{itemize}
\item
$x \mo^d \sqr A$    iff     $\ex U\ni x~ (U$ is open in $\X~ \&~ \fo
y\in U-\{ x\}~ y \mo^d  A)$;
\item
$x \mo^d \Di A$ iff     $\fo U\ni x~ (U$ is open in $\X~ \Ra~ \ex
y\in U-\{ x\}~ y \mo^d  A)$.
\end{itemize}
\end{lem}

\begin{defi}\label{lt1}
A \emph{local $T_1$-space} (or a {\em $T_D$-space} \cite{AT62}) is a topological space,
in  which every point is {\em locally closed}, i.e, closed in some 
neighbourhood. 
\end{defi}

Note that a point $x$ in an Alexandrov space $N(W,R)$ is closed iff it is minimal (i.e., $R^{-1}(x)=\{x\}$); $x$ is 
locally closed iff $R(x)\cap R^{-1}(x)=\set{x}$. Thus $N(F)$ is local
$T_1$ iff $F$ is a poset.

\begin{lem}\label{L310}\cite{E1}
For a topological space $\X$
\begin{enumerate}
\item
$\X \mo^d  \bbK4\mbox{  iff  } \X\mbox{ is local }T_1;$
\item
$\X \mo^d  \Di\top\mbox{  iff  } \X\mbox{ is dense-in-itself}.$
\end{enumerate}
\end{lem}

\begin{defi}
A Kripke frame (W,R)  is called {\em  weakly transitive} if
$R\circ R\sbq  \baR$.
\end{defi}
It is obvious that the weak transitivity of $R$ is equivalent to the transitivity of $\baR$.

\begin{propo}\label{P32}\cite{E1}~
(1) $(W,R) \mo  \bbK4^\circ\mbox{  iff  } (W,R)\mbox{ is weakly transitive};$\\
(2) $\bbK4^\circ$ is Kripke-complete.
\end{propo}

\begin{lem}\label{L33}\cite{E1}~
(1) Let $F = (W,R)$ be a Kripke $\bbS4$-frame, and let $R^\circ := R-I_W$,  $F^\circ := (W,R^\circ)$.
Then ${\bf Ld}(N(F)) = \BL(F^\circ)$.\\
(2)
Let $F=(W,R)$ be a weakly transitive irreflexive Kripke frame, and let $\baF=: (W, \baR)$ be its
reflexive closure. Then $ \BL\Bd (N(\baF))=\BL(F) $.\\
(3)
If $\bLa=\BL(\cC)$, for some class $\cC$ of weakly transitive irreflexive Kripke frames, then $\bLa$
is d-complete.
\end{lem}
\begin{proof}
(1) Note that $R^\circ(x)$ is the smallest punctured neighbourhood of $x$ in the space $N(F)$. So the inductive d-truth definition in a topological model $(N(F),\vp)$ coincides with the inductive truthdefinition in the Kripke model $(W,R^\circ,\vp)$.

(2) Readily follows from (1), since $\baR$ is transitive and $(\baR)^\circ=R$.

(3)  Follows from (2). {\hspace*{\fill} }
\end{proof}

\begin{defi}
For a 1-modal formula $A$ we define
$A^\sharp$ as the formula obtained by replacing every occurrence of every subformula $\sqr B$ with $\overline{\sqr}B:=\sqr B\we B$.
For a 1-modal logic $\bLa$ its {\em reflexive fragment} is $^\sharp\bLa:=\{A\mid\bLa\vd A^\sharp\}$.
\end{defi}

\begin{propo}\label{cldl}\cite{BEG}~
(1) If $\bLa$ is a $\bbK4^\circ$-logic, then $^\sharp\bLa$ is an $\bbS4$-logic.
\\
(2) For any topological space $X$, ${\bf Lc}(\X)= \,^\sharp{\bf Ld}(\X)$,
\\
(3) For any weakly transitive Kripke frame $F$, ${\bf L}(\overline{F})=\, ^\sharp{\bf L}(F)$.
\end{propo}

\begin{proof}
(1) It is clear that for a weakly transitive $\bLa$, $\overline{\sqr}$ satisfies the axioms of $\bbS4$, so  $^\sharp\bLa$ contains these axioms. Since $^\sharp$ distributes over implication, it follows that $^\sharp\bLa$ is closed under Modus Ponens. For the substitution closedness, note that for any variable $p$ and formulas $A,B$ $([B/p]A)^\sharp=[B^\sharp/p] A^\sharp$; thus $A\in\, ^\sharp\bLa$ implies $[ B/p]A\in\, ^\sharp\bLa$.Finally, since $(\sqr A)^\sharp= \overline{\sqr}\,A^\sharp$, it is clear that $A\in\, ^\sharp\bLa$   only if $\sqr A\in\,^\sharp\bLa$.

(2) By definitions,
$${\bf Lc}(X)\vd A \miff CA(X)\mo A,$$
$$^\sharp{\bf Ld}(X)\vd A \miff {\bf Ld}(X)\vd A^\sharp\miff DA(X)\mo A^\sharp.$$
Let us show that that $ CA(X)\nmo A$ iff $ DA(X)\nmo A^\sharp$. In fact, consider a topological model  $(X,\vp)$. We claim that
\[
\vp_c(B)=\vp_d(B^\sharp)\eqno(*)
\]
for any formula $B$. This is easily checked by induction, the crucial case is when $B=\sqr B_1$; then by definitions and the induction hypothesis we have:
\[
\vp_c(B)=\BI\vp_c(B_1)=\BI\vp_d(B_1^\sharp)= \boxdot\vp_d(B_1^\sharp)\cap \vp_d(B_1^\sharp)=\vp_d(\overline{\sqr}\,B_1^\sharp)= \vp_d(B^\sharp).
\]
The claim (*) implies that $\vp_c(A)\neq X$ iff   $\vp_d(A^\sharp)\neq X$ as required.

(3) On the one hand,
$${\bf L}(\overline{F})=\BL(MA(\overline{F}))= {\bf L} (CA(N(\overline{F}))= {\bf Lc}(N(\overline{F})).$$
On the other hand, by Lemma \ref{L33}(2),
\[
{\bf L}(F)= {\bf Ld}(N(\overline{F})),
\]
and we can apply (2) to $N(F)$.{\hspace*{\fill} }
\end{proof}

Let us give some examples of d-complete logics.
\begin{enumerate}
\renewcommand{\labelenumi}{(\arabic{enumi})}
\item ${\bf Ld}(\mbox{all ~topological ~spaces})=
\bbK4^\circ$. This was proved by L. Esakia in the 1970s and published in \cite{E1}.
\item
${\bf Ld}(\mbox{all ~local~}T_1\mbox{-spaces})= \bbK4$. This is also a result from \cite{E1}.

\item
${\bf Ld}(\mbox{all }T_0\mbox{-spaces})= \bbK4^\circ + p \land \romb (q \land \romb p) \to \romb p \lor \romb(q \land \romb q)$. This result is from \cite{BEG2011}.

\item
L. Esakia \cite{E} also proved that G\"odel - L\"ob logic
${\bf GL} := \bK + \sqr (\sqr p\sp p)\sp  \sqr p$
is the derivational logic of the class of all topological scattered spaces (a space is {\em
scattered} if each its nonempty subset has an isolated point).
\item
The papers \cite{A87}, \cite{A88}, \cite{Bl} give a complete description of
d-logics of
ordinals with the interval topology: ${\bf Ld}(\al)$ is
either ${\bf
GL}$ (if $\al\geq\om^\om$),  or
${\bf  GL}+ \sqr^n \bot$  (if $\om^{n-1}\leq\al<\om^n$). In particular, 
${\bf Ver} := \BK + \sqr \bot$ is the d-logic of any finite ordinal (and of any discrete space).
\item
The well-known ``difference logic" 
\cite{S80}, \cite{DR93}
${\bf DL} := \bbK4^\circ  +  \Di \sqr p \sp  p,$
is determined by Kripke frames with the difference relation:
${\bf DL} = \BL(\{(W,\neq_W)\mid W\neq\emp\}),$
where  $\neq_W:=W^2-I_W$; hence by \ref{cldl}, ${\bf DL}$ is the d-logic of the class of all
trivial topological spaces. However, for any particular trivial space $\X$, ${\bf Ld}(\X)\neq {\bf
DL}$. Moreover, ${\bf Ld}(\X)$  is not finitely axiomatizable for any infinite trivial $\X$\cite{KudShap}; this surprising result is easily proved by a standard technique using Jankov formulas (cf. \cite{Kudinov08}).
\item
In \cite{Sh00} it was proved that
${\bf Ld}\mbox{(all 0-dimensional separable metric spaces)}=\bbK4.$
All these spaces are embeddable in $\BR$ \cite{K66}.
\item
In \cite{Sh00} it was also proved that for any dense-in itself separable metric space $\X$, ${\bf
Ld}(\X)=\bbD4$; this was a generalization of an earlier proof \cite{Sh90} for $\X=\BQ$.  A more elegant proof for $\BQ$ is in \cite{LB1}.
\item
Every extension of $\bbK4$ by a set of closed axioms is a d-logic of some subspace of $\BQ$ \cite{BLB}. This gives us a
continuum of d-logics of countable metric spaces.

\item
In \cite{Sh90} ${\bf Ld}(\BR^2)$ was axiomatized and it was also proved that the d-logics of $\BR^n$ for $n\geq 2$ coincide. 
We will simplify and extend that proof in the present chapter.
\item
$\mathbf{Ld}(\BR)$ was described in \cite{Sh00}; for a simpler completeness proof cf. \cite{LB2}.
\item $\mathbf{Ld} (\mbox{all Stone spaces})= \bbK4$ and $\mathbf{Ld} (\mbox{all weakly scattered Stone spaces})= \bbK4 + \romb \top \to \romb \Box \perp$,
cf. \cite{BEG2010}.

\item
d-logics of special types of spaces were studied in \cite{BEG}, \cite{LB1}. They include submaximal, perfectly disconnected, maximal, weakly scattered and some others.
\end{enumerate}

However, not all extensions of $\bbK4^\circ$ are d-complete. In fact, the formula $p\sp \Di p$ never
can be d-valid, because  $\Bd Y = \emp$  for a singleton $Y$. So every  extension of $\bbS4$ is
d-incomplete, and thus Kripke completeness does not imply d-completeness.

\begin{propo}\label{P42}
Let $F=(\om^*,\prec)$ be the ``standard irreflexive transitive tree", where  $\om^*$  is the set of all finite sequences in
$\om$; $\al\prec\be$ iff  $\al$ is a proper initial segment of  $\be$. Then
$$\bD4 = \BL(F)=\BL\Bd (N(\overline{F})) = \BL\Bd (\cD),$$
where  $\cD$ denotes the class of all dense-in-themselves local $T_1$-spaces.
\end{propo}
\begin{proof}
The first equality is well known \cite{VB83}; the second one holds by \ref{L33}. By \ref{L310},
$\bbD4$ is d-valid  exactly in spaces from $\cD$. So $N(\overline{F})\in\cD$, $\bbD4\sbq\BL\Bd (\cD)
$, and the third equality follows.{\hspace*{\fill} }
\end{proof}

\section{Adding the universal modality and the difference modality}
\label{sec:basics}
\renewcommand{\labelenumi}{(\arabic{enumi})}

Recall that the {\em universal modality} $[\fo]$ and the {\em difference modality} $[\neq]$ correspond to Kripke frames with the universal and the difference relation.
So (under a valuation in a set $W$)
these modalities are interpreted in the standard way:
\begin{align*}
x\mo[\fo]A &\miff~ \fo y\in W~y\mo A;~ &x\mo[\neq]A &\miff \fo y\in W~(y\neq x \Ra y\mo A).
\end{align*}

The corresponding dual modalities are denoted by $\langle\ex\rangle$ and $\DE$.

\begin{defi}
For a $[\fo]$-modal formula $A$
we define the $[\neq]$-modal formula $A^u$
by induction:
\begin{equation*}
A^u:=A\mbox{ for }A\mbox{ atomic},\ \
 (A\sp B)^u:=A^u\sp B^u,\ \
([\fo]B)^u :=[\neq]B^u\we B^u.
\end{equation*}
\end{defi}

We can consider 2-modal topological logics obtained from ${\bf Lc}(\X)$ or ${\bf Ld}(\X)$ by adding
the universal or the difference modality\footnote{So we extend the definitions of the d-truth or the
c-truth by adding the item for $[\fo]$ or $[\neq]$.}. Thus for a topological space $\X$ we obtain
four 2-modal logics : ${\bf Lc}_\fo(\X)$ (the {\em closure universal (cu-) logic}), ${\bf
Ld}_\fo(\X)$ (the {\em derivational universal (du-) logic}), ${\bf Lc}_{\neq}(\X)$ (the {\em closure
difference (cd-) logic}), ${\bf Ld}_{\neq}(\X)$ (the {\em derivational difference (dd-) logic}).
Similar notations (${\bf Lc}_\fo(\cC)$ etc.) are used for logics of a class of spaces
$\cC$, and respectively we can define four kinds of topological completeness (cu-, du-, cd-, dd-) for
2-modal logics.

cd-logics were first studied in \cite{Gab01},
 cu-logics in \cite{Sh99}, du-logics in
\cite{LB2}, but dd-logics have never been addressed so far.

For a $\sqr$-modal logic $\BL$ we define the 2-modal logics
\begin{align*}
\BL\BD &:=\BL*\BD\BL+[\neq]p\we p\ri\sqr p,
~
&\BL\BD^+ &:=\BL*\BD\BL+[\neq]p\ri\sqr p,
\\
\BL\BU&:=\BL*\bS5+[\fo]p\ri\sqr p.
\end{align*}
Here we suppose
that $\bS5$ is formulated in the language with $[\fo]$ and $\BD\BL$ in the language with $[\neq]$.
The following is checked easily:
\begin{lem}\label{L41}
For any topological space $\X$,
\[
{\bf Lc}_\fo(\X)\spq \bbS4\BU,\quad {\bf Ld}_\fo(\X)\spq\bbK4^\circ\BU,\quad {\bf Lc}_{\neq}(\X)\spq
\bbS4\BD,\quad {\bf Ld}_{\neq}(\X) \spq\bbK4^\circ\BD^+.
\]
\end{lem}

\begin{defi}
For a 1-modal Kripke frame $F=(W,R)$ we define 2-modal frames 
$F_\fo:=(F,W^2),~ F_{\neq}:=(F,\neq_W)$
and modal logics $\BL_\fo(F):=\BL(F_\fo),~ \BL_{\neq}(F):=\BL(F_{\neq})$.
\end{defi}

Sahlqvist theorem \cite{CZ97} implies
\begin{propo}~
The logics $\bbS4\BU,~ \bbK4^\circ\BU, ~ \bbS4\BD,~\bbK4^\circ\BD^+ $ are Kripke complete.
\end{propo}
Using the first-order equivalents of the modal axioms for these logics (in particular, Proposition
\ref{P32}) we obtain

\begin{lem}\label{L46}~ For a rooted Kripke frame $G=(W,R,S)$
\begin{enumerate}
\item
$
G\mo\bbS4\BU\miff R\mbox{ is a quasi-order }\&~S=W^2,
$
\item
$
G\mo\bbK4^\circ\BU \miff R\mbox{ is weakly transitive }\&~S=W^2,
$
\item
$ G\mo\bbS4\BD \miff R\mbox{ is a quasi-order }\&~\overline{S}=W^2, $
\item
$ G\mo\bbK4^\circ\BD^+\miff R\mbox{ is weakly transitive }\&~\overline{S}=W^2\&~R\sbq S. $
\end{enumerate}
Also note that $\overline{S}=W^2\miff\neq_W\sbq S$.
\end{lem}

\begin{defi}\label{F0}
A rooted Kripke $\lgc{K4^\circ D}^+$-frame described by
Lemma \ref{L46} (4) is called {\em basic}. The class of these frames is denoted by $\mathfrak{F}_0$.
\end{defi}

Next, we easily obtain the 2-modal analogue to Lemma \ref{L33}.

\begin{lem}\label{L47}~
(1)
Let $F$ be an $\bbS4$-frame. Then
\[
{\bf Ld}_{\neq}(N(F))=\BL_{\neq}(F^\circ),~ {\bf Ld}_\fo(N(F))=\BL_\fo(F^\circ).
\]
(2)
Let $F$ be a weakly transitive irreflexive Kripke frame. Then
\[
{\bf Ld}_{\neq}(N(\overline{F}))=\BL_{\neq}(F),~ {\bf Ld}_\fo(N(\overline{F}))=\BL_\fo(F).
\]
(3)
Let $\cC$ be a class of weakly transitive irreflexive Kripke 1-frames. Then $\BL_{\neq}(\cC)$
is dd-complete, $\BL_{\fo}(\cC)$
is du-complete.
\end{lem}
Let us extend the translations $(-)^\sharp$, $(-)^u$ to 2-modal formulas.

\begin{defi}
$(-)^u$ translates $(\sqr, [\fo])$-modal formulas to $(\sqr, [\neq])$-modal formulas so that $([\fo]B)^u =[\neq]B^u\we B^u$
and $(-)^u$ distributes over the other connectives.

Similarly, $(-)^\sharp$ translates $(\sqr,[\neq])$-modal formulas and
$(\sqr,[\fo])$-modal formulas to formulas of the same kind, so that 
$(\sqr\,B)^\sharp =\sqr\,B^\sharp\we B^\sharp$
and $(-)^\sharp$ distributes over the other connectives.
\begin{align*}
^u\bLa&:=\{A\mid A^u\in\bLa\}\mbox{ for a }(\sqr,[\fo])\mbox{-modal logic } \bLa\mbox{ (the {\em
universal fragment}),}\\
^\sharp\bLa&:=\{A\mid A^\sharp\in\bLa\}\mbox{ for a }(\sqr,[\neq])\mbox{-  or a }
(\sqr,[\fo])\mbox{-modal }\bLa\mbox{ (the {\em reflexive fragment}),}\\
\,^\sharp\phantom{}^u\bLa&:= \,^\sharp(^u\bLa)\mbox{ for a }(\sqr,[\neq])\mbox{-modal  } \bLa\mbox{
(the {\em reflexive universal fragment}).}
\end{align*}
\end{defi}

\begin{propo}\label{P49}~
(1)
The map $\bLa\mt\,^\sharp\bLa$ sends $\lgc{K4^\circ D^+}$-logics to $\lgc{S4U}$-logics.\\
(2)
The map $\bLa\mt\,^u\bLa$ sends $\lgc{K4^\circ D^+}$-logics to
$\lgc{ K4^\circ U}$-logics and $\lgc{S4D}$-logics to $\lgc{S4U}$-logics.
\\(3)
The map $\bLa\mt\,^\sharp\phantom{}^u\bLa$ sends $\lgc{K4^\circ D^+}$-logics to
$\lgc{S4U}$-logics.
\\(4)
For a topological space $\X$
\[
\BL\Bc_{\neq}(\X)=\,^\sharp\BL\Bd_{\neq}(\X),~ \BL\Bd_{\fo}(\X)=\,^u\BL\Bd_{\neq}(\X),~ \BL\Bc_{\fo}(\X)=\,^u\BL\Bc_{\neq}(\X)={} ^\sharp\BL\Bd_{\fo}(\X).
\]
\\(5)
For a weakly transitive Kripke frame $F$
\[
\BL_{\neq}(\overline{F})=\,^\sharp\BL_{\neq}(F),~
~\BL_{\fo}(F)=\,^u\BL_{\neq}(F),
~\BL_{\fo}(\overline{F})=\,^\sharp\BL_{\fo}(F)=
\,^\sharp\,^u\BL_{\neq}(F).
\]
\end{propo}

Proposition \ref{P49}\,(4) implies that dd-logics are the most
expressive of all kinds of the logics we consider.

\begin{coro}
If $\BL\Bd_{\neq}(\X)= \BL\Bd_{\neq}(\YY)$ for spaces $\X,\YY$,  then all the other
logics (du-, cu-, cd-, d-, c-) of these spaces coincide.
\end{coro}


Let
\[
AT_1:=  \DA p \ri \DA \sqr p, \quad
AC:=  \UA(\sqr p \lor \sqr\neg p) \ri \UA p \lor \UA \lnot p.
\]

\begin{propo}\label{P412}
For a topological space $\X$
 \begin{enumerate}
\renewcommand{\labelenumi}{(\arabic{enumi})}
\item $\X \models^d \Di\top$ iff $\X$ is dense-in-itself;
\item $\X \models^d AT_1$ iff $\X \models^c AT_1$ iff $\X$ is a $T_1$-space;
\item $\X\mo^d AC^{\sharp}$ iff $\X \models^c AC$ iff $\X$ is connected.
 \end{enumerate}
\end{propo}

\begin{proof}
 (1) and the first equivalence in (2) are trivial. The first equivalence in (3) follows from
\ref{P49}(4). The remaining ones are checked easily, cf. \cite{kudinov_2006},
\cite{Sh99}.{\hspace*{\fill} }
\end{proof}

For a $\sqr$-modal logic $\lgc{L}$ put
\[
 \begin{array}{rcl}
\lgc{LD^+T_1} := \lgc{LD^+} + AT_1, &\qquad&
\lgc{LD^+T_1C} := \lgc{LD^+} + AT_1 + AC^{\sharp u}.
 \end{array}
\]
Also put
\[
\lgc{KT_1}:=\lgc{K4D^+T_1},~\lgc{DT_1}:=\lgc{D4D^+T_1},
~\lgc{DT_1C}:=\lgc{D4D^+T_1C}.
\]

\begin{propo}\cite{kudinov_2006}\label{pr:DS_AT1}
If $F = (W, R, R_D)$ is basic, then \mbox{$F \models AT_1$} iff all $R_D$-irreflexive points are strictly $R$-minimal iff $R_D \circ R \subseteq R_D$.
\end{propo}

%

\begin{rem}\rm
Density-in-itself is expressible in cd-logic and dd-logic by the formula $DS:= \DA p \impl \romb p$,
So for any space $\X$,
$\X\mo^c DS\miff\X\mo^d DS\miff \X\mo^d\Di\top$.
It is known that $DS$ axiomatizes dense-in-themselves spaces in cd-logic \cite{kudinov_2006}.
However, in dd-logic this axiom is insufficient:
${\bf Ld}_{\neq}$(all dense-in-themselves spaces) $= \lgc{D4^\circ D^+}=\lgc{K4^\circ D^+}+\Di\top$, and it is
{\em stronger} than $\lgc{K4^\circ D^+}+DS$. (To see the latter, consider a
singleton Kripke frame, which is $R_D$-reflexive, but $R$-irreflexive.) Therefore
$\lgc{K4^\circ D^+}+DS$ is dd-incomplete.
\end{rem}

\begin{rem}\rm
Every $T_1$-space is a local $T_1$-space, so the dd-logic of all  $T_1$-spaces contains $\sqr
p\ri\sqr\sqr p$. However, $\lgc{K4^{\circ}D^+T_1}\,\nvd\,\sqr p\ri\sqr\sqr p$.
In fact, consider a 2-point frame $F:=(W,\neq_W,W^2)$. It is clear that $F\mo\lgc{K4^{\circ}D^+}$.
Also $F\mo AT_1$, by Proposition \ref{pr:DS_AT1}, but $F\nmo\sqr p\ri\sqr\sqr p$, since  $\neq_W$ is
not transitive.

It follows that $\lgc{K4^{\circ}D^+T_1}$ is dd-incomplete; $T_1$-spaces are actually axiomatized by ${\bf KT_1}$ (Corollary \ref{C713}).
\end{rem}





\newcommand{\Lcu}{\mathbf{Lc_\forall}}
\newcommand{\Lcd}{\mathbf{Lc_{\ne}}}
\newcommand{\Ldu}{\mathbf{Ld_\forall}}
\newcommand{\Ldd}{\mathbf{Ld_{\ne}}}

Let us give some examples of du-, cu- and cd-complete logics.
\begin{enumerate}
\item $\Lcu(\hbox{all spaces}) = \lgc{S4U}$.
 \item $\Lcu(\hbox{all connected spaces}) = \Lcu(\Real^n) = \lgc{S4U} + AC$  for any $n\ge 1$  \cite{Sh99}\footnote{The paper \cite{Sh99} contains a stronger claim: $\Lcu(\X)=\lgc{S4U} + AC$ for any connected dense-in-itself separable metric $\X$. However, recently we found a gap 
in the proof of Lemma 17 from that paper. Now we state the main result only for the case ${\X}=\Real^n$; a proof can be obtained by applying the methods of the present Chapter, but we are planning to publish it separately.}  

 \item $\Ldu(\hbox{all spaces}) = \lgc{S4D}$ \cite{DR93}.
 \item $\Lcd(\X) = \lgc{S4DT_1+DS}$, where $\X$ is a 0-dimensional separable metric space \cite{kudinov_2006}.
 \item $\Lcd(\BR^n)$ for any $n\ge 2$ is finitely
axiomatized in \cite{kudinov_2005};  all these logics coincide.
 \item $\Ldu(\BR)$ is finitely axiomatized in \cite{LB2}.
\end{enumerate}

\section{dd-completeness of $\lgc{K4^\circ D^+}$ and some of its extensions}
This section contains some simple arguments showing that there are many dd-complete bimodal logics.

All formulas and logics in this section are $(\sqr,\DA)$-modal. An arbitrary Kripke frame for
$(\sqr,\DA)$-formulas is often denoted by $(W,R,R_D)$.

\begin{lem}\label{L35}~
(1)
Every weakly transitive Kripke 1-frame is a p-morphic image of  some irreflexive weakly transitive Kripke 1-frame. \\
(2)
Every rooted $\lgc{K4^\circ D^+}$-frame is a p-morphic image of some $R$- and $R_D$-irreflexive rooted
$\lgc{K4^\circ D^+}$-frame.
\end{lem}
\begin{proof}
(1) Cf. \cite{E1}.

(2) Similar to the proof of (1).
For $F=(W,R,R_D)\in \mathfrak{F}_0$ put
\[
W_r:=\{a\mid aR_Da\},~ W_i=W-W_r, \quad \tilde{W}:= W_i\cup (W_r \ti \{ 0,1\}).
\]
Then we define the relation $\tilde{R}$  on $\tilde{W}$ such that 
\begin{align*}
 (b,j) \tilde{R}a      &\mbox{ ~~iff ~~ }   bRa,&
a \tilde{R} (b,j)      &\mbox{ ~~iff ~~ }   aRb,\\
(b,j) \tilde{R} (b',k)      &\mbox{ ~~iff ~~ }   bRb'\ \&\  b\ne b' \vee  b= b'\ \&\  j\ne k,&
a\tilde{R}  a'      &\mbox{ ~~iff ~~ }   aRa'.
\end{align*}

Here  $a,a'\in W_i;~ b,b'\in W_r; ~j,k\in\{ 0,1\}$. So we duplicate all $R_D$-reflexive points making them irreflexive (under
both relations). 
It follows that $\tF:=(\tW,\tilde{R},\neq_{\tW})\in \mathfrak{F}_0$ and $\tilde{R}$ is
irreflexive; the map $f:\tW \to W$ sending $(b,j)$ to $b$ and $a$ to itself (for $b\in W_r, ~ a\in
W_i$) is a p-morphism $\tF\dri F$. {\hspace*{\fill} }
\end{proof}

\begin{propo}\label{P53}
Let $\Gamma$ be a set of closed 2-modal formulas, $\bLa:= \lgc{K4^\circ D^+}+\Gamma$. Then
\begin{enumerate}
\item
$\bLa$ is Kripke complete.
\item
$\bLa$ is dd-complete.
\end{enumerate}
\end{propo}

\begin{proof}
(1) $\lgc{K4^\circ D^+}$ is axiomatized by Sahlqvist formulas. One can easily check that (in the
minimal modal logic) every closed formula is equivalent to a positive formula, so
we can apply Sahlqvist theorem.

(2) Suppose $A\nin\bLa$. By (1) and the Generation lemma there exists a rooted Kripke 2-frame $F$
such that $F\mo L$ and $F\nmo A$. Then by Lemma \ref{L35}, for some irreflexive weakly transitive
1-frame $G=(W,R)$ there is a p-morphism $(G,\neq_W)\dri F$. By the p-morphism lemma 
$(G,\neq_W)\nmo A$ and $(G,\neq_W)\mo\bLa$ (since $\Gamma$ consists of closed formulas). Hence by Lemma
\ref{L47}, $\bLa\sbq\lgc{Ld}_{\neq}(N(\overline{G}))$, $A\nin\lgc{Ld}_{\neq}(N(\overline{G}))$. 
{\hspace*{\fill} }
\end{proof}

\begin{rem}\rm
Using Proposition \ref{P53} and the construction from \cite{BLB} one can prove that there is a continuum of dd-complete logics. Such a claim is rather weak, because Proposition \ref{P53} deals only with Alexandrov spaces. In section 7 we will show how to construct many dd-complete logics of metric spaces.
\end{rem}

\section{d-morphisms and dd-morphisms; extended McKinsey - Tarski's Lemma}
In this section we recall the notion of a d-morphism (a 
validity-preserving map for d-logics) and introduce dd-morphisms, the analogues of d-morphisms for dd-logics. This is the main technical tool in the present chapter. Two basic lemmas are proved here, an analogue of 
McKinsey--Tarski's lemma on dissectability for d-morphisms and the Glueing lemma.

The original McKinsey--Tarski's lemma \cite{MT} states the existence of a c-morphism (cf. Remark \ref{R61} ) from an arbitrary separable dense-in-itself metric space onto a certain quasi-tree of depth 2. The 
separability condition is actually redundant \cite[Ch.~3]{RS} (note that the latter proof is quite different from  \cite{MT}\footnote{\label{foot:7}Recently P. Kremer \cite{Kremer} has showed that $\bbS4$ is {\em strongly complete} w.r.t. any dense-in-itself metric space. His proof uses much of the construction from \cite{RS}.}). But c-morphisms preserve validity only for c-logics, and unfortunately, the constructions by 
McKinsey--Tarski and Rasiowa--Sikorski cannot be used for d-morphisms. So we need another construction to prove a stronger form of McKinsey--Tarski's lemma. 

\begin{defi}\label{D51}
Let $\X$  be a topological space, $F=(W,R)$ a transitive Kripke frame. A map $f: X\lora W$ is called a {\em
d-morphism} from  $\X$ to $F$ if $f$ is open and continuous as a map $\X\lora N(\overline{F})$ and also satisfies
\begin{align*}
\mbox{r-density}:&~~\fo w\in W( wRw\Ra f^{-1} (w)\sbq \Bd f^{-1} (w)),\\
\mbox{i-discreteness}:&~~\fo w\in W( \neg wRw\Ra f^{-1} (w)\cap \Bd f^{-1} (w)=\emp).
\end{align*}
If $f$ is surjective, we write $f: \X\dri^d F$.
\end{defi}

\begin{propo}\label{P52}\cite{BEG}~
(1)
$f$ is a d-morphism from $\X$ to $F$ iff $2^f$ is a homomorphism from $MA(F)$ to $DA(\X)$.\\
(2)
If $f:  \X\dri^d F$, then $\BL\Bd (\X) \sbq  \BL(F)$.
\end{propo}

\begin{coro}\cite{Sh90}\label{L52}
A map $f$ from a  topological space $\X$ to a finite transitive Kripke frame $F$ is a d-morphism iff
$$\fo w\in W~ \Bd f^{-1} (w) = f^{-1}(R^{-1}(w)).$$
\end{coro}

\begin{proof}
$2^f$ preserves Boolean operations. It is a homomorphism of modal algebras iff it preserves diamonds, i.e., iff for any $V\sbq W$,
$$f^{-1}(R^{-1}(V))=\Bd f^{-1}(V).$$
Inverse images and $\Bd$ distribute over finite unions, so the above equality holds for any (finite) $V$ iff it holds for singletons, i.e.,
\[\hspace{3.8cm}\ f^{-1}(R^{-1}(w))= \Bd f^{-1}(w).\qedhere \]
\end{proof}

\begin{rem}\label{R61}\rm
For a space $\X$ and a Kripke $\bbS4$-frame $F=(W,R)$ one can also define a {\em  c-morphism} $\X\lora F$ just as an open and continuous map $f: \X\lora N(F)$. So every d-morphism to an $\bbS4$-frame is a c-morphism. It is well known \cite{RS} that $f:  X\lora W$ is a c-morphism iff $2^f$
is a homomorphism $MA(F) \lora CA(\X)$.
Again for a finite $F$ this is equivalent to
$$
\fo w\in W~\BC f^{-1}(w) = f^{-1}(R^{-1}(w)). 
$$
\end{rem}

\begin{lem}\label{Restlem}
If $f:\X\dri^d F$ for a finite frame $F$ and $\Y\sbq\X$ is an open subspace, then $f|Y$ is a d-morphism.
\end{lem}

\begin{proof}
We apply Proposition \ref{P52}. Note that $f|Y$ is the composition $f \cd j$, where
$j:Y\hookrightarrow X$ is the inclusion map. Then $2^{ f|Y}=2^j\cd 2^f$. Since $2^f$ is a
homomorphism $MA(F)\lora DA(\X)$, it remains to show that $2^j$ is a homomorphism $DA(\X)\lora
DA(\Y)$, i.e., it preserves the derivation:
$j^{-1}(\Bd V)=\Bd_Yj ^{-1}(V)$, or
$\Bd V\cap Y=\Bd_Y(V\cap Y)$,
which follows from \ref{dY}. 
{\hspace*{\fill} }
\end{proof}

\begin{defi}\label{D53}
A set $\gamma$ of subsets of a topological space
$\X$ is called \emph{dense} at $x\in X$ if every neighbourhood of $x$
contains a member of $\gamma$.
\end{defi}

\begin{figure}[h]
 \centering
 \includegraphics[width=0.4\textwidth]{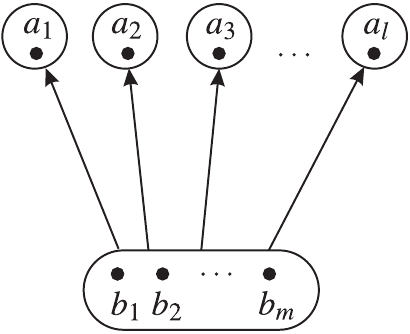}
 \caption{Frame $\Phi_{ml}$.}
 \label{fig:Phiml}
\end{figure}

\begin{propo}\label{P54}
For $m>0,~ l>0$ let
$\Phi_{ml}$ be a ``quasi-tree" of height 2, with singleton maximal clusters and
An
 $m$-element
root cluster (Fig. \ref{fig:Phiml}). For $l=0, ~m>0, ~ \Phi_{ml}$ denotes an $m$-element cluster.


Let $\X$ be a dense-in-itself separable metric space, $B\ssb X$ a closed nowhere dense set.
Then there exists a d-morphism $g: \X\dri^d \Phi_{ml}$ with the following properties:
\begin{enumerate}
\item
$B \sbq  g^{-1}(b_1)$;
\item
every $g^{-1}(a_i)$ (for $i\leq l$ )  is a union of a set $\al_i$  of disjoint
open balls, which  is
dense at any point of $g^{-1}(\{ b_1,...,b_m\}$).
\end{enumerate}
\end{propo}
\begin{proof}
Let $X_1,\dots, X_n,\dots$ be a countable base of $\X$ consisting of open balls. We construct sets
$A_{ik}, ~B_{jk}$  for $1\leq i\leq l, ~1\leq j\leq m, ~k\in\om$, with the following properties:
\begin{enumerate}
\renewcommand{\labelenumi}{(\arabic{enumi})}
\item
$A_{ik}$ is the union of a finite set $\al_{ik}$ of nonempty open balls whose closures are disjoint;
\item
$\BC A_{ik} \cap  \BC A_{i'k} = \emp$   for  $i\not= i'$ ;
\item
$\al_{ik} \sbq  \al_{i,k+1};~A_{ik} \sbq  A_{i,k+1}$;
\item
$B_{jk}$  is finite;
\item
$B_{jk} \sbq  B_{j,k+1};$
\item
$A_{ik} \cap  B_{jk} = \emp ;$
\item
$X_{k+1} \sbq\bigcup\limits_{i=1}^l A_{ik}\Ra \al_{i,k+1} = \al_{ik}, ~B_{j,k+1}
=
B_{jk};$
\item
if $X_{k+1}\not\sbq \bigcup\limits_{i=1}^l A_{ik}$, there are closed
nontrivial balls
$P_1,\dots,P_l$ such that for any $i$, $j$
$$P_i \sbq  X_{k+1}-A_{ik},~ 
\al_{i,k+1} = \al_{ik} \cup  \{ \BI P_i\},~
(B_{j,k+1} - B_{jk}) \cap  X_{k+1}\not= \emp; $$
\item
$A_{ik} \sbq  X-B$;
\item
$B_{jk} \sbq  X-B$;
\item
$j\not= j'\Ra   B_{j'k} \cap  B_{jk} = \emp$ .
\end{enumerate}
We carry out both the construction and the proof by induction on $k$.

Let $k=0$. $(X-B)$  is infinite, since 
it is nonempty and open in a
dense-in-itself
$\X$.
Take distinct points $v_1,\dots,v_l \nin B$ and disjoint closed
nontrivial
balls
$ Z_1,\dots,Z_l \ssb X-B$ with centres at $v_1,\dots,v_l$ respectively (see Fig.\ref{fig:McT-constraction}).

\begin{figure}[h]
 \centering
 \includegraphics[width=0.4\textwidth]{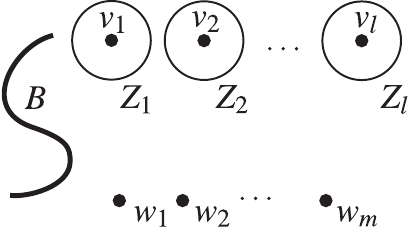}
 \caption{Case k = 0}
 \label{fig:McT-constraction}
\end{figure}

Put
 $$\al_{i0}: =\{ \BI Z_i\}; ~A_{i0}: = \BI Z_i;$$
then $Z_i = \BC A_{i0}$. As above, since $(X-B) -\bigcup\limits^l_{i=1}Z_i$  is nonempty and open, it
is infinite. Pick distinct $w_1,\dots,w_m\in X-B$ and
put $B_{j0}: = \{ w_j\}$. Then the required properties hold for
$k=0$.

At the induction step we construct
$A_{i,k+1}, B_{j,k+1}$.
Put $Y_k :=\bigcup\limits ^l_{i=1}A_{ik}$ and consider two cases.\\
(a)  $X_{k+1}\sbq Y_k$. Then
put:
$$\al_{i,k+1}: = \al_{ik};~ A_{i,k+1}: = A_{ik}; ~B_{j,k+1}: = B_{jk}.$$
(b)   $X_{k+1}\not\sbq  Y_k$. Then
$X_{k+1} \not\sbq \BC Y_k$.
In fact, $X_{k+1} \sbq \BC Y_k$ implies $X_{k+1} \sbq  \BI \BC Y_k=Y_k$, since $X_{k+1}$ is open and by (1) and (2).
So we put
$$W_0 := X_{k+1} - \BC Y_k -\bigcup^m_{j=1}B_{jk},~
W := W_0 - B.$$
Since $(X_{k+1} - \BC Y_k )$ is nonempty and open and every $B_{jk}$ is finite by
(4), $W_0$  is also open and nonempty (by the density of $\X$).  By the assumption of
\ref{P54}, $B$ is closed, and thus $W$ is open.

$W$ is also nonempty. In fact, otherwise $W_0 \sbq  B$, and then $W_0 \sbq  \BI B=\emp$  (since
$B$ is nowhere dense by the assumption of \ref{P54}).

Now we argue similarly to the case $k=0$.
Take disjoint closed nontrivial balls $P_1,\dots,P_l\ssb W$. Then $W - \bigcup\limits^l_{i=1}P_i$
is infinite, so we pick distinct $b_{1,k+1},\dots,
b_{m,k+1}$ in this set and put
$$B_{j,k+1}: = B_{jk} \cup  \{ b_{j,k+1}\},~
\al_{i,k+1}: = \al_{ik} \cup  \{ \BI P_i\},~ A_{i,k+1}: = A_{ik} \cup  \BI P_i.$$
In the case (a) all the required properties hold for $(k+1)$ by the construction.

In the case (b) we have to check only (1), (2), (6), (8)--(11).

(8) holds, since by construction we have
\begin{align*}
P_i\ssb W\ssb& X_{k+1}-\BC Y_k\ssb X_{k+1}-A_{ik};\\
b_{j,k+1}\in W \sbq&  X_{k+1},~ b_{j,k+1}\in (B_{j,k+1}-B_{jk}).
\end{align*}

(1): From IH it is clear that $\al_{i,k+1}$ is a finite set of open balls and
their closures are disjoint; note 
that  $P_i \cap  \BC A_{ik} = \emp$, since
$P_i \sbq  W  \sbq  -\BC A_{ik}.$

(2): We have
\begin{align*}
&\BC A_{i,k+1} \cap  \BC A_{i^\prime,k+1} = (\BC A_{ik} \cup  P_i) \cap  (\BC
A_{i^\prime k}\cup  P_{i^\prime}) =\\
=&(\BC A_{ik} \cap  \BC A_{i^\prime k}) \cup  (\BC A_{ik} \cap  P_{i'}) \cup  (\BC
A_{i^\prime k} \cap  P_i) \cup  (P_i \cap  P_{i'})
 =\BC A_{ik} \cap  \BC A _{i^\prime k} = \emp
\end{align*}
by IH and by the construction; note that
 $P_i, P_i'  \sbq  W \sbq  -\BC Y_k$.

(6): We have
$$
A_{i,k+1} \cap  B_{j,k+1} = (A_{ik} \cap  B_{jk} ) \cup  (\BI P_i \cap  
\set{b_{j,k+1}}) \cup  ( A_{ik} \cap  \set{ b_{j,k+1}}) \cup  (\BI P_i \cap  B_{jk}) = \emp
$$ 
by IH and since $b_{j,k+1}\nin P_i,~ b_{j,k+1}\in W \sbq X-Y_k$, $P_i \ssb W \sbq
X-B_{jk}$ .

(9): We have
$A_{i,k+1} = A_{ik} \cup  \BI P_i \sbq  -B,$ 
since $A_{ik} \sbq  -B$ by IH, and $P_i \ssb W \sbq  -B$ by the construction.

Likewise, (10)  follows from $B_{jk} \sbq  -B$ and $b_{j,k+1}\in W \sbq  -B$.

To check (11), assume $j \not=  j'$. We have
$B_{j',k+1} \cap  B_{j,k+1} = B_{j'k} \cap  B_{jk}$,
since
$b_{j^\prime,k+1} \not=  b_{j,k+1} , ~ b_{j,k+1}\in W \sbq  -B_{j^\prime k}$
and
$b_{j^\prime,k+1}\in W\sbq -B_{jk}$. Then apply IH.

Therefore the required sets $A_{ik}, B_{jk}$ are constructed. Now put
$$\al_i: =\bigcup_k \al_{ik},~
A_i:= \bigcup \al_i=\bigcup_k A_{ik}, ~
B_j: =\bigcup_k B_{jk},$$
$$ B^\prime _1: = X - (\bigcup_i A_i \cup \bigcup_j B_j),$$
and define a map
$g: X\lora \Phi_{ml}$  as follows:
\[
g(x):=\left\{ \begin{array}{ll}
a_i & \mbox{if } x\in A_i, \\
b_j & \mbox{if } x\in B_j,~j\neq 1,\\
b_1 & \mbox{otherwise (i.e., for } x\in B^\prime _1). \\
\end{array} \right.
 \]
By (2), (3), (5), (6), (11),  $g$ is well defined; by (9), (10) $B \sbq  g^{-1}(b_1)$.

To prove that $g$  is a d-morphism, we check some other properties.
$$\leqno(12)\quad    X - \bigcup\limits^l_{i=1}A_i \sbq  \Bd B_j.$$

In fact, take an arbitrary $x\nin \bigcup\limits^l_{i=1}A_i  $ and show that $x\in \Bd B_j$, i.\/e., 
$$\leqno(13)\quad(U-\{ x\}) \cap  B_j \not=  \emp . $$
for any neighbourhood $U$  of $x$. First assume that $x\nin B_j$. Take a basic open $X_{k+1}$ such that $x\in X_{k+1} \sbq  U$. Then
\mbox{$X_{k+1}\not\sbq \bigcup\limits^l_{i=1}A_i$,} 
and (8) implies
$B_{j,k+1} \cap  X_{k+1} \not=  \emp .$
Thus
$B_j \cap  U \not=  \emp$.
So we obtain (13).

Suppose $x\in B_j$; then $x\in B_{jk}$  for some $k$. Since $\X$ is dense-in-itself and
$\set{X_1,\,X_2,\dots}$ is its open base, $\setdef[X_{s+1}]{s\geq k}$ is also an open base 
(note that every ball in $\X$ contains a smaller ball). 
So $x\in X_{s+1}\sbq U$ for some $s\ge k$. Since $x\nin \bigcup\limits^l_{i=1}A_i $, we have
$X_{s+1}\not\sbq\bigcup\limits^l_{i=1}A_i$,
and so $(B_{j,s+1}-B_{js}) \cap  X_{s+1} \not=  \emp$ by (8); thus $(B_j-B_{js}) \cap  U \not=  \emp$.
Now $x\in B_{jk} \sbq  B_{js}$ implies (13).
$$\leqno(14)\qquad    \Bd B_j \sbq  X -\bigcup\limits^l_{i=1}A_i.$$
In fact,
$B_j\sbq -A_i$, by (3), (5), (6). So $\Bd B_j \sbq  \Bd (-A_i) \sbq  -A_i$, since $A_i $  is open.

Similarly we obtain 
$$\leqno(15)\qquad    \Bd B^\prime _1\sbq  X -\bigcup\limits^l_{i=1}A_i, \qquad    \Bd A_i \sbq  X- \bigcup\limits_{r\not=i}A_r.$$
Also note that 
$$\leqno(16)\qquad    A_i \sbq  \Bd A_i,$$
since $A_i$  is open, $\X$ is dense-in-itself. Similary to (12) we have
$$\leqno(17)\qquad    \al_i \mbox{ is dense at every point of  }B_j, B'_1\ \mbox{  (and thus $B_j,\ B'_1 \sbq  \Bd A_i$}).$$


To conclude that $g$ is a d-morphism, note that
$$g^{-1}(a_i) = A_i, ~g^{-1}(b_j) = B_j~ (\mbox{for }j\not= 1),~ g^{-1}(b_1) =
B^\prime _1,$$
and so by (15), (16), (17)
\begin{align*}
\Bd g^{-1}(a_i) &= \Bd A_i = X - \bigcup\limits_{r\not=i}A_r = g^{-1}(R^{-
1}(a_i)),
\end{align*}
and by (12), (14), (15)
\begin{align*}
\Bd g^{-1}(b_j) &= \Bd B_j = X - \bigcup\limits^l_{i=1}A_i = g^{-1}(R^{-
1}(b_j))\ \ \mbox{(for $j\ne 1$)},\\ \qquad\qquad\qquad
\Bd g^{-1}(b_1) &= \Bd B^\prime _1= X - \bigcup\limits^l_{i=1}A_i = g^{-1}(R^{-
1}(b_1)). \qedhere
\end{align*}
\end{proof}

For the proof see Appendix.

\begin{lem}\label{refroot}
Assume that
\begin{enumerate}
\item
$\X$ is a dense-in-itself separable metric space,
\item
$B\subset X$ is closed nowhere dense,
\item
$F = C\cup F_1\cup \dots\cup F_l$ is a $\lgc{D4}$-frame, where $C = \{ b_1,\dots,b_m\} $ is a non-degenerate
root cluster,
$F_1,\dots,F_l$  are the subframes generated by the successors of $C$,
\item
for any nonempty open ball $U$ in $\X$, for any $i\in\{1,\ldots,l\}$ there exists a d-morphism $f_i^U:U\dri^d F_i$.
\end{enumerate}
Then there exists $f:\X\dri^d F$ such that $f(B)=\{b_1\}$.
\end{lem}
\begin{proof}
First, we construct
$g: \X\dri^d \Phi_{ml}$ according to Proposition \ref{P54}. Then $B\sbq g^{-1}(b_1)$ and
$A_i = g^{-1}(a_i)$ is the union of a set $\al_i$ of disjoint open balls.
Then put
\begin{equation}\label{Eq:fPhiml}
f(x):=\left\{ \begin{array}{ll}
g(x) & \mbox{if } g(x)\in C, \\
f_i^U(x)& \mbox{if } x\in U,~U\in\al_i.\\
\end{array} \right.
\end{equation}

Since $g$ and $f_i^U$ are surjective, the same holds for $f$. So let us show 
$$\Bd f^{-1}(a) = f^{-1}(R^{-1}(a))$$
($R$ is the accessibility relation on $F$).
First suppose $a\in C$. Then (since  $g$ is a d-morphism)
$$\Bd f^{-1}(a) = \Bd g^{-1}(a) = g^{-1}(C) = f^{-1}(C) = f^{-1}(R^{-1}(a)). $$
Now suppose $a \notin C$, $I = \setdef[i]{a \in F_i}$, and let $R_i$ be the accessibility relation on $F_i$. We have:
\begin{align*}
f^{-1}(a) =&\bigcup\limits_{i\in I}\bigcup\limits_{U\in\al_i}(f_i^U)^{-1}(a),&
R^{-1}(a) =& C \cup  \bigcup\limits_{i\in I} R_i^{-1}(a),
\end{align*}
and so 
$$ f^{-1}(R^{-1}(a)) =
     g^{-1}(C) \cup \bigcup\limits_{i\in I}\bigcup\limits_{U\in\al_i} (f_i^U)^{-1}(R_i^{-1}(a)).
     $$     
Since $f_i^U$ is a d-morphism,
\begin{equation}\label{eq:L69:2}
f^{-1}(R^{-1}(a)) =
g^{-1}(C) \cup\bigcup\limits_{i\in I}\bigcup\limits_{U\in\al_i} \Bd_U((f_i^U)^{-1}(a))
 \sbq  g^{-1}(C) \cup  \Bd f^{-1}(a).
\end{equation}

Let us show that
\begin{equation}\label{eq:L69:3} g^{-1}(C) \sbq  \Bd f^{-1}(a).
 \end{equation}

In fact, let $x\in g^{-1}(C)$. Since 
$\al_i$ is dense at
$x$, every  neighbourhood of $x$
contains some $U\in\al_i$. Since $f^U_i$ is
surjective,
$f(u) = f^U_i (u) = a$ for some $u\in U$. Therefore,  $x\in\Bd f^{-1}(a)$. 

(\ref{eq:L69:2}) and (\ref{eq:L69:3}) imply $f^{-1}(R^{-1}(a)) \sbq  \Bd f^{-1}(a)$.
Let us prove the converse:
\begin{equation}\label{eq:L69:5}  \Bd f^{-1}(a) \sbq  f^{-1}(R^{-1}(a)).
\end{equation}

We have $A_j \cap  f^{-1}(a)=\emp$  for $j\notin I$  and $A_j$   is open,
hence
$A_j \cap  \Bd f^{-1}(a)=\emp$.
Thus  $\Bd f^{-1}(a) \sbq  g^{-1}(C) \cup  A_i$.
Now $g^{-1}(C)\sbq f^{-1}(R^{-1}(a))$ by (\ref{eq:L69:2}), so it remains to show that for any $i\in I$
\begin{equation}\label{eq:L69:7}
\Bd f^{-1}(a)\cap A_i\sbq f^{-1}(R^{-1}(a)).
\end{equation}

To check this, consider any $x\in\Bd f^{-1}(a)\cap A_i$. Then $x\in U$ for some $U\in \al_i$, and thus by \ref{dY} and  (\ref{eq:L69:2}) 
$x\in\Bd f^{-1}(a)\cap U=\Bd_U(f^{-1}(a)\cap U)=\Bd_U(f^U_i)^{-1}(a)\sbq f^{-1}(R^{-1}(a))$. 
This implies (\ref{eq:L69:7}) and completes the proof of (\ref{eq:L69:5}).
{\hspace*{\fill} }
\end{proof}

Recall that $\partial$ denotes the boundary of a set in a topological space: $\partial A := \BC A - \BI A$.

\begin{lem}\label{Glue}{\em (Glueing lemma)}
Let $\X$ be a local $T_1$-space satisfying 

\noindent
(a) $X=X_1\cup Y\cup X_2$ for closed nonempty subsets $X_1,Y,X_2$ such that
\begin{itemize}
\item
$X_1\cap X_2=X_1\cap\BI Y= X_2\cap\BI Y =\emp$,
\item
$\partial X_1\cup\partial X_2=\partial Y$,
\item
$\Bd\BI Y=Y$ (i.e., $Y$ is regular and dense in-itself).
\end{itemize}
or

(b) $X=X_1\cup X_2$ is a nontrivial closed partition.

Let $F=(W,R)$ be a finite $\bbK4$-frame, 
$F_1=(W_1,R_1),~ F_2=(W_2,R_2)$ its generated subframes such that $W=W_1\cup W_2$ and suppose there are
d-morphisms $f_i:\X_i\dri^d F_i$, $i=1,\,2$, where $\X_i$ is the subspace of $\X$ corresponding to $X_i$.

In the case (a) we also assume that $F_1, F_2$ have a common maximal cluster $C$, $f_i(\partial
X_i)\sbq R^{-1}(C)$ for $i=1,2$ and there is $g:\BI Y\dri^d C$ (where $C$ is regarded as
a frame with the universal relation, $\BI Y$ as a subspace of $\X$). Then
$f_1\cup f_2\cup g:\X\dri^d F$ in the case (a), $f_1\cup f_2:\X\dri^d F$ in the case
(b).\footnote{$f_1\cup f_2$ is the map $f$ such that $f|X_i=f_i$; similarly for $f_1\cup f_2\cup g$.}
\end{lem}

\begin{figure}[h]
 \centering
 \includegraphics[width=0.5\textwidth]{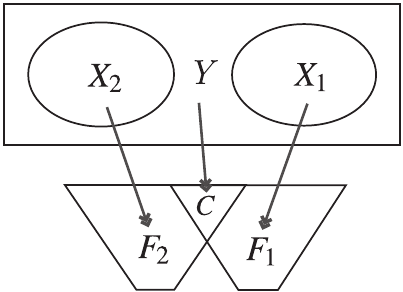}
 \caption{Case (a)}
 \label{fig:glueing_lemma}
\end{figure}

\begin{proof}
Let $f:= f_1\cup f_2\cup g $ (or $f:= f_1\cup f_2$), $ F_i=(W_i,R_i)$,
$\Bd:=\Bd_X,~\Bd_i:=\Bd_{X_i}$.
For $w\in W$ there are four options.

(1) $w\in W_1-W_2$. Then
$\Bd f^{-1}(w)=\Bd f_1^{-1}(w)=\Bd_1 f_1^{-1}(w) = f_1^{-1}(R_1^{-1}(w))$ (since $X_1$ is closed and $f_1$ is a d-morphism).
It remains to note that $R_1^{-1}(w)=R^{-1}(w)\sbq W_1-W_2$, and thus
$f_1^{-1}(R_1^{-1}(w))= f^{-1}(R^{-1}(w)).$

(2) $w\in W_2-W_1$. Similar to the case (1).

(3) $w\in (W_1\cap W_2)-C$ in the case (a) or $w\in W_1\cap W_2$ in the case (b).  Then $f^{-1}(w)=f_1^{-1}(w)\cup f_2^{-1}(w)$, so similarly to (1), 
\[
\Bd f^{-1}(w)=  \Bd_1 f_1^{-1}(w)\cup \Bd_2 f_2^{-1}(w) = f_1^{-1}(R_1^{-1}(w))\cup f_2^{-1}(R_2^{-1}(w))= f^{-1}(R^{-1}(w)).
\]

(4) $w\in C$ in case (a). First note that $\Bd g^{-1}(w)=Y$. In fact,
$g$ is a d-morphism onto the cluster $C$, so
$\Bd_{\BI Y}  g^{-1}(w)= g^{-1}(C)=\BI Y$. Hence 
$\BI Y\sbq \Bd  g^{-1}(w )\sbq \Bd\BI Y=Y$, and thus
$$Y=\Bd\BI Y\sbq\Bd\Bd g^{-1}(w ) \sbq \Bd  g^{-1}(w )$$
by \ref{L310}(2).
Next, since $X_1$, $X_2$ are closed and $f_1$, $f_2$ are d-morphisms we have
\begin{align*}
\Bd f^{-1}(w) &= \Bd f_1^{-1}(w)\cup \Bd f_2^{-1}(w)\cup \Bd g^{-1}(w) 
=
\Bd_1 f_1^{-1}(w)\cup \Bd_2 f_2^{-1}(w)\cup Y =&\\
&= f_1^{-1}(R_1^{-1}(w))\cup f_2^{-1}(R_2^{-1}(w))\cup Y= f^{-1}(R^{-1}(w)). &\qedhere
\end{align*}
\end{proof}

The case (b) of the previous lemma can be generalized as follows.
\begin{lem}\label{L67}
Suppose a topological space $\X$  is the disjoint union of open subspaces:
$\X=\bigsqcup\limits_{i\in I}\X_i$. Suppose a Kripke $\bbK4$-frame  $F$ is the
union of its generated subframes: $F=\bigcup\limits_{i\in I}F_i$
and suppose $f_i:\X_i\dri^{d}F_i$. Then $\bigcup\limits_{i\in I}f_i:\X\dri^{d}F$.
\end{lem}


\begin{defi}\label{def_pmorphism}
 Let $\X$ be a topological space,
 $F= (W, R, R_D)$ be a frame. Then a surjective map $f: X \lora W$ is called a
 {\em dd-morphism}
 (in symbols, $f:\X \dri^{dd} F$) if
\begin{enumerate}
 \item $f:\X \dpmor (W, R)$ is a d-morphism ;
 \item $f:(X, \ne_X) \pmor (W, R_D)$ is a p-morphism of Kripke frames.
\end{enumerate}
\end{defi}

\begin{lem}\label{lem_pmorphism}
 If $f:\X \dri^{dd} F$, then $\BL\Bd_{\neq}(\X) \subseteq \BL(F)$ and for any closed  2-modal $A$
\[
\X\mo A\Lra F\mo A.
\]
\end{lem}

\begin{proof}
Similar to \ref{P52} and \ref{p1}.
{\hspace*{\fill}
}
\end{proof}

\begin{defi}
A set-theoretic map $f:X\lora Y$ is called {\em n-fold at} $y\in Y$ if $|f^{-1}(y)|=n$;\footnote{$|\ldots|$ denotes the cardinality.}
$f$ is called {\em manifold at} $y$ if it $n$-fold for some $n>1$.
\end{defi}

\begin{propo}\label{lem:dmortodDmor}~
(1) Let $G = (X, \ne_X)$,  $F = (W, S)$ be Kripke frames such that $\overline{S}=W^2$, and
let $f : X \lora W$ be a surjective function. Then
$$f : G \pmor F \quad\miff\quad f\mbox{ is manifold exactly at }S\mbox{-reflexive points of }F.$$
(2) 
Let $\X$ be a $T_1$-space, $F= (W, R, R_D)$ a rooted
$\lgc{KT_1}$-frame, $f:\X\dri^d(W,R)$. Then $f:\X\dri^{dd}F$ iff for any strictly $R$-minimal $v$
\[
v R_D v\Lra f\mbox{ is manifold at }v.
\]
(3) 
If $\X$ is a $T_1$-space,  $f: \X \dpmor F=(W, R)$ and $R^{-1}(w)\neq\emp$ for any $w\in W$, then $f:\X \dri^{dd} F_{\fo}$,
where $F_{\fo}:= (W, R, W^2)$.
\end{propo}

\begin{proof}
(1) Note that
$f$ is a p-morphism iff for any $x\in X$
\[
f(X-\{x\})=S(f(x))=
\begin{cases}
W & {\rm if~}f(x)Sf(x),\nonumber\\
W-\{f(x)\} & {\rm otherwise.}\nonumber
 \end{cases}
\]

(2) By (1), $f:\X\dri^{dd}F$ iff
\[
\fo v\in W(v R_D v\Lra \abs{f^{-1}(v)} > 1).
\]
The latter equivalence holds whenever $ R^{-1}(v)\neq\emp$. In fact, then by Corollary \ref{L52},
$\Bd f^{-1} (v) = f^{-1}(R^{-1}(v)) \neq\emp$, and thus $f^{-1}(v)$ is not a singleton (since $\X$ is
a $T_1$-space). $ R^{-1}(v)\neq\emp$ also implies $vR_Dv$, by Proposition \ref{pr:DS_AT1}.

(3) follows from (2).
{\hspace*{\fill}
}
\end{proof}

After we have proved the main technical results, in the next sections we will study dd-logics of specific spaces.




\section{$\lgc{D4}$ and $\lgc{DT_1}$ as logics of zero-dimensional dense-in-themselves spaces}

In this section we will prove the d-completeness of $\lgc{D4}$ and 
dd-completeness of $\lgc{DT_1}$ w.r.t. zero-dimensional spaces. The proof follows rather easily from the previous section and an additional technical fact (Proposition \ref{P61}) similar to the McKinsey--Tarski lemma.

Recall that a (nonempty) topological space $\X$ is called
{\em zero-dimensional}  
if clopen sets
constitute its open base \cite{A77}.
Zero-dimensional $T_1$-spaces with a countable base are subspaces of the Cantor
discontinuum, or of
the
set of irrationals  \cite{K66}.

\begin{lem}\label{L71}
Let $\X$ be a zero-dimensional dense-in-itself Hausdorff space.
Then for any $n$ there exists a nontrivial open
partition $\X=\X_1\sqcup\ldots\sqcup\X_n$, in which every $\X_i$ is also a zero-dimensional dense-in-itself Hausdorff space.
\end{lem}

\begin{proof}
It is sufficient to prove the claim for $n=2$ and then apply induction.
A dense-in-itself space cannot be a singleton, so there are two different points $x,y\in X$. Since $\X$ is $T_1$ and zero-dimensional, there exists a clopen $U$ such that $x\in U,~y\nin U$. So
$X=U\cup (X-U)$ is a nontrivial open partition.
The Hausdorff property, density-in-itself, zero-dimensionality are inherited for open subspaces.
{\hspace*{\fill}}
\end{proof}

\begin{propo}\label{P61}
Let $\X$ be a zero-dimensional dense-in-itself metric space, $y\in X$. Let
$\Psi_l$ be the
frame consisting of an irreflexive root $b$ and its reflexive successors
$a_0,\dots,a_{l-1}$ (Fig. \ref{fig:Phiml_irreflexive_root}).

\begin{figure}[h]
 \centering
 \includegraphics[width=0.4\textwidth]{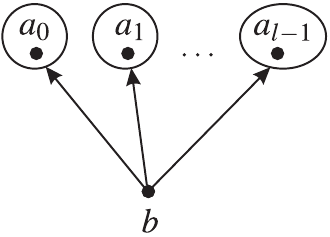}
 \caption{Frame $\Psi_l$.}
 \label{fig:Phiml_irreflexive_root}
\end{figure}

Then there exists $f: \X\dri^d \Psi_l$  such that $f(y) = b$  and for every $i$ 
there is an open partition of $f^{-1}(a_i)$, which is dense at y.
\end{propo}
\begin{proof}
Let $O(a,r) := \setdef[x\in X]{\rho(a,x)<r}$, where $\rho$ is the distance in $\X$.

There exist clopen sets $Y_0, Y_1, \dots$ such that
$$\set{y}\ssb\dots\ssb Y_{n+1} \ssb Y_n \ssb\dots Y_1 \ssb Y_0 = X$$ \
and
$Y_n \sbq  O(y, 1/n)$ for $n >0$.

These $Y_n$  can be easily constructed by induction. Then
$$\bigcap\limits_n Y_n =\set{y},
\ \hbox{and}\
X-\{ y\} =\bigsqcup\limits_n X_n, $$
where $X_n = Y_n-Y_{n+1}$. Note that the $X_n$ are nonempty and
open, $X_n \sbq  O(y, 1/n)$ for $n>0$.

Now define a map $f: X\lora\Psi_l$ as follows:
\[
f(x)=\left\{ \begin{array}{ll}
a_{r(n)} & \mbox{if } x\in X_n; \\
b & \mbox{if } x=y,\\
\end{array} \right.
 \]
where $r(n)$ is the remainder of dividing $n$  by $l$; it is clear that $f$ is surjective.

Let us show that for any $x$,
\[
x\in \Bd f^{-1}(u) \mbox{   iff     }f(x)Ru.
\eqno(*)
\]

(i)  Assume that $u=a_j$. Then
 $f^{-1}(u)=\bigcup\limits_n X_{nl+j}$,
and
$$f(x)Ru \mbox{ iff     }(f(x)=b \mbox{ or }  f(x)=u).$$
To prove `if' in (*), consider two cases.

1. Suppose $f(x)=u, ~x\in X_{nl+j}$. Since $X_{nl+j}$ is nonempty and open, it is
dense-in-itself, and thus $x\in \Bd X_{nl+j}\sbq \Bd f^{-1}(u)$.

2. Suppose $f(x)=b$, i.e. $x=y$. Then  $x\in \Bd f^{-1}(u)$, since
$X_{nl+j}\sbq O(y, 1/n)$.

The previous argument also shows that
$\{ X_{nl+j}\mid n\geq 0\}$ is an open partition of $f^{-1}(a_j)$, which
is dense at $y$.

To prove `only if', suppose $f(x)Ru$ is not true. Then
$f(x) = a_k$ for some $k \not=  j$,  and so for some $n$, $x\in
X_n$,  $X_n\cap f^{-1}(u)=\emp$. Since $X_n$ is open, $x\nin \Bd f^{-1}(u)$.

(ii) Assume that $u = b$. Then $f^{-1}(u) = \{ y\}$, and so $\Bd f^{-1}(u) =
\emp  = f^{-
1}(R^{-1}(u))$.
{\hspace*{\fill} }\end{proof}

\begin{propo}\label{P64}
Let $\X$ be a  zero-dimensional dense-in-itself separable metric space, $F$ a
finite rooted $\bD4$-frame.
Then there exists a d-morphism $\X\dri^dF$, which is 1-fold at the root of $F$
if this root is irreflexive.
\end{propo}
\begin{proof}
By induction on the size of $F$.

(i) If $F$ is a finite cluster, the claim follows from Proposition \ref{P54}.

(ii) If $F = C\cup F_1\cup \dots\cup F_l$, where $C = \{ b_1,\dots,b_m\} $ is a non-degenerate
root cluster,
$F_1,\dots,F_l$  are the subframes generated by the successors of $C$, we can apply Lemma \ref{refroot}. In fact, every open ball $U$ in $\X$ is zero-dimensional and dense-in-itself.

(iii)  Suppose $F=\breve{b}\cup F_0\cup\dots\cup F_{l-1}$, where $b$ is an irreflexive
root of $F$, $F_i$ are the subframes generated by the successors of $b$.
There exists
$g:  X\dri^d \Psi_l$  by \ref{P61}, with an arbitrary
$y\in X$. Then $g^{-1}(a_i)$ is a union of a set $\al_i$ of disjoint open
sets, and $\al_i$  is dense
at
$y$. If $U\in\al_i$, then by IH, there exists
$f^U_i: U\dri^d F_i$. Put
\[
f(x)=\left\{ \begin{array}{ll}
b & \mbox{ if } x=y; \\
f_i^U(x)& \mbox{ if } x\in U,~U\in\al_i.\\
\end{array} \right.
 \]
Then similarly to Lemma \ref{refroot} it follows that $f: X\dri^d F$.

Finally note that if the root of $F$ is irreflexive, the first step of the construction is case (iii),
so the preimage of the root is a singleton.
{\hspace*{\fill} }
\end{proof}
\begin{theorem}\label{T65}
If $\X$ is a zero-dimensional dense-in-itself separable metric space, then $\BL\Bd (\X)=\bD4$.
\end{theorem}
\begin{proof}
By Propositions \ref{P64}
and \ref{P52} $\BL\Bd (\X) \sbq \BL (F)$ for any finite rooted $\bD4$-frame $F$, thus $\BL\Bd (\X) \sbq  \bD4$, since $\bD4$ has the fmp. By Lemma \ref{L310}
$\bD4\sbq\BL\Bd (\X)$ . {\hspace*{\fill} }\end{proof}

\begin{lem}\label{P741}
Let $\X$ be a zero-dimensional dense-in-itself separable metric
space, $F$ a finite $\bD4$-frame. Then there exists a d-morphism $\X\dri^d F$, which is 1-fold at all strictly minimal points.
\end{lem}

\begin{proof}
 $F = F_1\cup\ldots\cup F_n$ for different finite rooted $\bD4$-frames $F_i$.
By Lemma \ref{L71}, $\X=\X_1\sqcup\ldots\sqcup\X_n$ for  zero-dimensional dense-in-themselves subspaces
$\X_i$, which are also metric and separable.
By Proposition \ref{P64}, we construct $f_i:\X_i\dri^d F_i$. Then by Lemma \ref{L67}, $\bigcup\limits_{i=1}^nf_i:\X\dri^d F$. Every strictly minimal point of $F$ is an irreflexive root of a unique $F_i$, so its preimage is a singleton.
{\hspace*{\fill} }
\end{proof}

\begin{propo}\label{P:dd-morph_zero_dem}
Let $\X$ be a zero-dimensional dense-in-itself separable metric space, $F \in
\mathfrak{F}_0$ a finite $\lgc{DT_1}$-frame. Then there exists a dd-morphism $\X\dri^{dd}F$.
\end{propo}
\begin{proof}
We slightly modify the proof of the previous lemma.
Let $F=(W,R,R_D)$, $G=(W,R)$. Then $G = G_1\cup\ldots\cup G_n$ for different cones $G_i$.
We call $G_i$ {\em special} if its root is strictly $R$-minimal
and $R_D$-reflexive.
We may assume that exactly $G_1,\ldots,G_m$ are special.
Then we count them twice and present $G$ as
$G_1\cup G'_1\cup\ldots\cup G_m\cup G'_m\cup G_{m+1}\cup\ldots\cup G_n$, where $G'_i=G_i$ for $i\leq m$ (or as $G_1\cup G'_1\cup\ldots\cup G_m\cup G'_m$ if $m=n$).

Now we can argue as in the proof of Lemma \ref{P741}. By Lemma \ref{L71},
$\X=\X_1\sqcup\X'_1\sqcup\ldots\sqcup\X_m\sqcup \X'_m\sqcup \X_{m+1}\sqcup\ldots\sqcup \X_n$ for
zero-dimensional dense-in-itself separable metric $\X_i, \X'_i$. By Proposition \ref{P64}, we
construct the maps $f_i:\X_i\dri^d G_i$, $f'_i:\X'_i\dri^d G'_i$, which are 1-fold at irreflexive
roots; hence by Lemma \ref{L67}, $f:\X\dri^d G$ for $f:= \bigcup\limits_{i=1}^nf_i\cup
\bigcup\limits_{i=1}^mf'_i $.

Every strictly minimal point $a \in G$ is an irreflexive root of a unique $G_i$.
If $a$ is $R_D$-irreflexive, then $G_i$ is not special, so
$f^{-1}(a)=f_i^{-1}(a)$ is a singleton.
If $a$ is $R_D$-reflexive, then $G_i$ is special, so $f^{-1}(a)= f_i^{-1}(a)\cup (f'_i)^{-1}(a)$, and thus $f$ is 2-fold at $a$.
Therefore, $f:\X\dri^{dd} F$
by Proposition \ref{lem:dmortodDmor}.
{\hspace*{\fill} }
\end{proof}

\begin{lem}\label{LF1}
 Let $M = (W, R, R_D, \vp)$ be a rooted Kripke model over a basic frame\footnote{Basic frames were defined in Section 4.} validating
$ AT_1$, $\Psi$  a set of 2-modal formulas closed under subformulas. Let $M'=(W',R',R_D',\th')$ be a
filtration  of $M$ through $\Psi$ described in Lemma \ref{L26}\footnote{Recall that $R^\prime$ is
the transitive closure of  $\uR$, $R_D'=\underline{R_D}$.}. Then the frame $(W',R',R_D')$ is also
basic and validates $AT_1$.
\end{lem}

\begin{proof}
In fact, $R'$ is transitive by definition. For any two different $a, b\in W'$ we have $a R_D'b$,
since $xR_Dy$ for any $x\in a,~y\in b$ (as $F\in\mathfrak{F}_0$).

Next, note that if $a$ is  $ R_D'$-irreflexive, then $a=\{x\}$ for some $R_D$-irreflexive $x$. In this case,
since $(W, R, R_D)\mo AT_1$, there is no $y$ such that $yRx$ (Proposition \ref{pr:DS_AT1}),
hence $(R')^{-1}(a)=\emp$, and thus $(W',R',R_D')\mo AT_1$.

Finally, $R'\sbq R_D'$. In fact, all different points in $F'$ are $ R_D'$-related, so it remains to
show that every $ R_D'$-irreflexive point is $R'$-irreflexive. As noted above, such a point is a
singleton class $x^\sim=\{x\}$, where $x$ is $R_D$-irreflexive. Then $x$ is $R$-minimal, so in $W'$
there is no loop of the form $x^\sim\uR x_1\uR\ldots\uR x^\sim$, and thus  $x^\sim$ is
$R'$-irreflexive. {\hspace*{\fill} }
\end{proof}

By a standard argument Lemma \ref{LF1} implies

\begin{theorem}\label{Tfmp}
Every logic of the form $\lgc{KT_1}+A$, where $A$ is a closed
2-modal formula, has the finite model property.
\end{theorem}

 \begin{proof}
 Let $L$ be such a logic and suppose $L\,\nvd\, B$. By Proposition \ref{P53} $L$ is Kripke complete,
 so by the Generation lemma there is a rooted Kripke frame $F=(W,R,R_D)$ such that $F\mo L,~F\nmo B$.
 Then $F$ is basic by definition. Let $M=(F,\th)$ be a Kripke model over $F$ refuting $B$. Let $\Psi$
 be the set of all subformulas of $A$ or $B$, and let us construct the filtration
 $M'=(W',R',R_D',\th')$ of $M$ through $\Psi$ as in Lemmas \ref{L26}(2) and \ref{LF1}. By the previous
 lemma, $F':= (W',R',R_D')\mo\lgc{KT_1}$.

 By the Filtration lemma, $M'\nmo B$. By the same lemma, the truth of $A$ is preserved in $M'$, so $F'\mo A$, since $A$ is closed. Therefore, $F'\mo L$.
 \end{proof}

\begin{theorem}\label{thm:DLogic_zerospace}
Let $\X$ be a zero-dimensional dense-in-itself separable metric space. Then  $\BL\Bd_{\neq}
(\X)=\lgc{DT_1}$.
\end{theorem}
\begin{proof}
For any finite $\lgc{DT_1}$-frame $F$ we have $\BL\Bd_{\neq} (\X) \sbq  \BL (F)$ by Proposition
\ref{P:dd-morph_zero_dem} and Lemma \ref{lem_pmorphism}. By the previous theorem, $\lgc{DT_1}$ has
the fmp, so  $\BL\Bd_{\neq} (\X) \sbq \lgc{DT_1}$. Since $\X\mo^d \lgc{DT_1}$ (Proposition
\ref{P412}), it follows that $\BL\Bd_{\neq} (\X) = \lgc{DT_1}$. {\hspace*{\fill} }
\end{proof}



\begin{propo}{\cite[Lemma 3.1]{BLB}}\label{prop:QtoF}
Every countable\footnote{In this chapter, as well as in \cite{BLB}, `countable' means `of cardinality at most $\aleph_0$'.} rooted $\lgc{K4}$-frame  is a
d-morphic image of a subspace of $\BQ$.
\end{propo}

To apply this proposition to the language with the difference
modality, we need to examine the preimage of the root
for the constructed morphism.
Fortunately, in the proof of Proposition \ref{prop:QtoF} in \cite{BLB}
the preimage of a root $r$ is a singleton iff $r$ is
irreflexive.

\begin{lem}\label{prop:QtoF_1fold}
 Let $F$ be a countable $\lgc{K4}$-frame. Then there exists a
d-morphism from a subspace of $\BQ$ onto $F$, which is 1-fold at
all strictly minimal points.
\end{lem}
\begin{proof}
Similar to Lemma \ref{P741}. We can present $F$ as a countable union of different cones
$\bigcup\limits_{i\in I}F_i$ and $\BQ$ as a disjoint union $\bigsqcup\limits_{i\in I}\X_i$ of
spaces homeomorphic to $\BQ$. By Proposition \ref{prop:QtoF} (and the remark after it), for
each $i$ there exists $f_i: \Y_i \dpmor F_i$ for some subspace $\Y_i \subseteq \X_i$ such that $f_i$ is
1-fold at the root $r_i$ of $F_i$ if $r_i$ is irreflexive. Now by Lemma \ref{L67} $f:=
\bigcup\limits_{i\in I}f_i: \bigsqcup\limits_{i\in I}\Y_i \dri^{d}F$, and $f$ is 1-fold at all
strictly minimal points of $F$ (i.e., the irreflexive $r_i$) --- since every $r_i$ belongs only to
$F_i$, so $f^{-1}(r_i)=f_i^{-1}(r_i)$. {\hspace*{\fill} }
\end{proof}

\begin{propo}\label{prop:subsetQtoF_ddmorph}
 Let $F$ be a countable $\lgc{KT_1}$-frame. Then there exists a
dd-morphism from a subspace of $\BQ$ onto $F$.
\end{propo}

\begin{proof}
Similar to Proposition \ref{P:dd-morph_zero_dem}. If $F = (W, R, R_D)$, the frame $G=(W,R)$ is a
countable union of different cones. There are two types of cones: non-special $G_i~ (i\in I)$ and
special (with strictly $R$-minimal and $R_D$-reflexive roots) $H_j~ (j\in J)$:
\[
G = \bigcup\limits_{i \in I} G_i \cup \bigcup\limits_{j \in J} H_j.
\]
Then we duplicate all special cones
\[
G = \bigcup\limits_{i \in I} G_i \cup \bigcup\limits_{j \in J} H_j \cup \bigcup\limits_{j \in J} H'_j
\]
and as in the proof of \ref{prop:QtoF_1fold}, construct $f: \bigsqcup\limits_{i\in
I}\Y_i\sqcup\bigsqcup\limits_{j\in J}\Z_j\sqcup\bigsqcup\limits_{j\in J}\Z'_j \dri^{d}F$. This map is
1-fold exactly at all  $R_D$-irreflexive points, so it is a dd-morphism onto $F$. {\hspace*{\fill}
}
\end{proof}

\begin{coro}\label{C713}
 $\mathbf{Ld}_{\ne}(\mbox{all $T_1$-spaces}) = \lgc{KT_1}$.
\end{coro}

\begin{proof}
Note that $\lgc{KT_1}$ is complete w.r.t.~countable frames and every subspace of $\BQ$ is $T_1$.
{\hspace*{\fill} }
\end{proof}

\begin{propo}\label{P714}
 Let $\bLa = \lgc{KT_1} + \Gamma$ be a consistent logic, where
$\Gamma$ is a set of closed formulas. Then $\bLa$ is
dd-complete w.r.t.\/ subspaces of $\BQ$.
\end{propo}
\begin{proof}
Since every closed formula is canonical, $\bLa$ is Kripke complete. So for every formula $A
\notin \bLa$ there is a frame $F_A$ such that $F_A \models \bLa$ and $F_A \nvDash A$. By
Proposition \ref{prop:subsetQtoF_ddmorph}, there is a subspace $\X_A \subseteq \BQ$ and  $f_A: \X_A \cpmor F_A$. Then $\X_A\nmo A,~\X_A\mo\bLa$ by Lemma \ref{lem_pmorphism}.
Therefore $\BL\Bd_{\neq} (\mathcal K) = \bLa$ for
$\mathcal K := \setdef[\X_A]{A \notin \bLa}$.
{\hspace*{\fill} }
\end{proof}

\begin{rem}\rm
A logic of the form described in Proposition \ref{P714} is dd-complete w.r.t. a set of subspaces of $\BQ$. This set may be non-equivalent to a single subspace. For example,
there is no subspace $\X\sbq\BQ$ such that $\lgc{KT_1}=\BL\Bd_{\neq}(\X)
$. In fact, consider
\[
A:=[\neq]\sqr\bot\we\sqr\bot.
\]
Then $A$ is satisfiable in $\X$ iff $\X\mo^d A$ iff $\X$ is discrete. So $A$ is consistent in
$\lgc{KT_1}$. Now if $\lgc{KT_1}=\BL\Bd_{\neq}(\X)$, then $A$ must be satisfiable in $\X$, hence  
$\X\mo^d A$; but  $\lgc{KT_1}\not\vd A$, and so we have a contradiction. 
\end{rem}


\section{Connectedness}
Connectedness was the first example of a property expressible in 
cu-logic, but not in c-logic. The corresponding connectedness axiom from \cite{Sh99} will be essential for our further studies.  
In this section we show that it is weakly canonical, i.e., valid in weak canonical frames --- a fact not mentioned in \cite{Sh99}.

\begin{lem}\label{L80}\cite{Sh99}
A topological space $\X$ is connected iff $\X\mo^c AC$, where
\[
AC :=[\fo](\sqr p\vee\sqr\neg p)\to [\fo] p\vee[\fo]\neg p.
\]
\end{lem}

For the case of Alexandrov topology there is an equivalent definition of connectedness in relational terms.
\begin{defi}\label{D75}
For a transitive Kripke frame $F = (W,R)$ we define the
{\em  comparability} relation $R^\pm:=R\cup R^{-1}\cup I_W$.
$F$ is called {\em  connected}  if the transitive closure of $R^\pm$ is universal.
A subset $V\sbq W$ is called {\em  connected}  in $F$ if the frame $F|V$ is connected.

A 2-modal frame $(W,R,S)$ is called  \emph{(R)-connected} if $(W,R)$ is connected.
\end{defi}
Thus $F$ is connected iff every two points $x,y$  can be connected by a {\em  non-oriented path}
(which we call just a {\em path}), a sequence of points $x_0x_1\ldots x_n$ such that
$x=x_0 R^\pm x_1\ldots R^\pm x_n=y$. 

From \cite{Sh99} and Proposition \ref{P49} we obtain

\begin{lem}\label{L801}
(1) For an $\bbS4$-frame $F$, the associated space $N(F)$ is connected iff $F$ is connected.\\
(2) For a $\bbK4$-frame $F$, $F_{\fo}\models AC^{\sharp u}$ iff $F$ is connected.
\end{lem}

%

\begin{lem}\label{mcl}
Let $M=(W,R,R_D,\th)$ be a rooted generated submodel of m-weak canonical model for a modal logic
$\bLa\spq{\bf K4D^+}$. Then
\begin{enumerate}
\item Every $R$-cluster in $M$ is finite of cardinality at most $2^m$.
\item
$(W,R)$ has finitely many $R$-maximal clusters.
\item
For each $R$-maximal cluster $C$ in $M$ there exists an $m$-formula $\be(C)$ such that:
$$\fo x\in M~ (M,x \mo  \be(C)\Lra x\in\overline{R}^{-1} (C)).$$
\end{enumerate}
\end{lem}

The proof is similar to \cite[Section 8.6]{CZ97}.

\begin{lem}\label{lem_KripkeCompletAC}
Every rooted generated subframe of a weak canonical
frame for a logic
$\bLa\spq{\bf K4D^+} + AC^{\sharp u} $ is connected.
\end{lem}

\begin{proof}
Let $M$ be a weak canonical model for $\bLa$, $M_0$ its rooted generated submodel with the frame
$F=(W,R,R_D)$, and suppose $F$ is disconnected. Then there exists a nonempty proper clopen subset $V$
in the space $N(W,\overline{R})$. Let $\Delta$ be the set of all $R$-maximal clusters in $V$ and put
\[
 B := \bigvee\limits_{C \in \Delta} \be(C).
\]
Then $B$ defines $V$ in $M_0$, i.e., $V=\overline{R}^{-1}(\bigcup\Del)$.
In fact, $\bigcup\Del\sbq V$ implies $\overline{R}^{-1}(\bigcup\Del) \sbq V $, since $V$ is closed. The other way round, $V\sbq\overline{R}^{-1}(\bigcup\Del)$, since for any $v\in V$, $\overline{R}(v)$ contains an $R$-maximal cluster $C\in\Del$, and $\overline{R}(v)\sbq V$ as $V$ is open.


So $w \models B$ for any $w \in V$, and since $V$ is open, $w \models \oB B$. By the same reason,
$w\models \oB\lnot B$ for any $w \not \in V $ . Hence
\[
M_0 \models \UA (\oB B \lor \oB \lnot B).
\]
By Proposition \ref{cm} all substitution instances of $AC$ are true in $M_0$. So we have
\[
M_0 \models \UA (\oB B \lor \oB \lnot B) \ri \UA B \lor \UA \lnot B,
\]
and thus
\[
M_0 \models \UA B \lor \UA \lnot B.
\]
This contradicts the fact that $V$ is a nonempty proper subset of
$W$. {\hspace*{\fill} }
\end{proof}

In d-logic instead of connectedness we can express some its local versions; they will be considered in
the next section. 
%

\section{Kuratowski formula and local 1-componency}
In this section we briefly study Kuratowski formula distinguishing $\BR$ from $\BR^2$ in d-logic. Here the main proofs are similar to the previous section, so most of the details are left to the reader.

\begin{defi}
We define {\em
Kuratowski formula} as
\[
Ku:= \quad \square (\oB p \vee  \oB \neg p) \ri  \square p \vee  \square \neg p.
\]
\end{defi}

The spaces validating $Ku$ are characterized as follows \cite{LB2}.

\begin{lem}\label{lem:th32_from_LB2}
For a topological space $\X$, $\X\mo^d Ku$ iff

\smallskip
\noindent
for any $x\in X$ and any open neighbourhood $U$ of $x$, if $U-\{x\}$ is a disjoint union
$V_1\cup V_2$ of sets open in the subspace $U-\{x\}$, then there exists a neighbourhood\footnote{In
\cite{LB2} neighbourhoods are supposed open, but this does not matter here, since every neighbourhood
contains an open neighbourhood.} $V \subseteq U$ of $x$ such that $V-\{x\}  \sbq V_1$ or  $V-\{x\}
\sbq V_2$.
\end{lem}

\begin{defi}
A topological space $\X$ is called {\em locally connected} if every  neighbourhood of any point $x$
contains a connected neighbourhood of $x$. Similarly, $\X$ is called {\em locally 1-component} if
every punctured neighbourhood of any point $x$ contains a connected punctured neighbourhood of $x$.
\end{defi}

It is well known \cite{A77} that in a locally connected space every  neighbourhood $U$ of any point $x$
contains a connected {\em open} neighbourhood of $x$ (e.g. the connected component of $x$ in $\BI U$).

\begin{lem}\label{LK1}
If $\X$ is locally 1-component, then
$\X\mo^d Ku$.
\end{lem}

The proof is straightforward, and we leave it to the reader.

\begin{lem}\label{LK2}~
(1) Every space d-validating $Ku$ has the following \emph{non-splitting property}:

(NSP) If an open set $U$ is connected, $x\in U$ and  $U-\{x\}$ is open, then $U-\{x\}$ is connected.\\
(2) 
Suppose $\X$ is locally connected and local $T_1$. Then (NSP) holds in $\X$ iff $\X$ is locally 1-component iff $\X\mo^d Ku$.
\end{lem}

\begin{proof}
(1) We assume $\X\mo^d Ku$ and check (NSP). Suppose $U$ is open and connected, $U^\circ:=U-\{x\}$ is
open, and consider a partition $U^\circ=U_1\cup U_2$ for open $U_1,U_2$. By 
\ref{lem:th32_from_LB2} there exists an open $V\sbq U$ containing $x$ such that $V\sbq\{x\}\cup U_1$
or $V\sbq\{x\}\cup U_2$. Consider the first option (the second one is similar). We have a partition
\[
U=(\{x\}\cup U_1)\cup U_2,
\]
and $\{x\}\cup U_1= V\cup U_1$, so $\{x\}\cup U_1$ is open. Hence by connectedness, $U= \{x\}\cup U_1$, i.e.,
$U^\circ = U_1$. Therefore,
$U^\circ $ is connected.

(2) It suffices to show that (NSP) implies the local 1-componency. Consider $x\in X$ and its neighbourhood $U_1$. Since $\X$ is local
$T_1$, $U_1$ contains an open neighborhood $U_2$, in which $x$ is closed, i.e., $\BC\{x\}\cap
U_2=\{x\}$. By the local connectedness, $U_2$ contains a connected open neighbourhood $U_3$, and
again $\BC\{x\}\cap U_3=\{x\}$; thus $U_3-\{x\}$ is open. Eventually, $U_3-\{x\}$ is connected, by
(NSP). {\hspace*{\fill} }
\end{proof}

\begin{rem}\rm
The {\em
($n$-th)
 generalized Kuratowski formula} is the following formula in variables $p_0,\dots,p_n$
$$Ku_n := \square \bigvee\limits^n_{k=0}\oB Q_k \ri
\bigvee\limits^n_{k=0}\square\neg Q_k,$$
where $Q_k := p_k \we\bigwedge\limits_{j\not= k}\neg p_j$.

The formula $Ku_1$ is related to the equality found by Kuratowski \cite{K22}:
$$(*)\quad\Bd ((x \cap  \Bd (-x)) \cup  (-x \cap  \Bd x)) = \Bd x \cap  \Bd (-x),$$
which  holds in every algebra $DA(\BR^n)$ for
$n>1$, but
not in $DA(\BR)$. This equality corresponds to the modal formula
$$Ku^\prime:=\quad   \Di ((p \we  \Di \neg p) \vee  (\neg p \we  \Di p)) \lrarr\Di p
\we  \Di \neg p,$$
and one can show that $\bD4 + Ku^\prime = \bD4 + Ku_1=\bD4 + Ku.$
\end{rem}

\begin{rem}\rm
The class of spaces validating $Ku_n$ is described in \cite{LB2}. In particular, it is valid in all locally  $n$-component spaces defined as follows.

A neighbourhood $U$ of a point $x$ in a topological space is called
{\em  $n$-component
at $x$} if the punctured neighbourhood $U-\set{x}$ has at most $n$
connected components.
A topological space is called {\em  locally  $n$-component} if
the $n$-component neighbourhoods at each of its point constitute a local base (i.e., every neighbourhood contains an $n$-component neighbourhood).
\end{rem}

\begin{lem}\label{L77}\cite{LB2}
For a transitive Kripke frame $(W,R)$

$(W,R) \mo  Ku$  iff  for any $R$-irreflexive $x$, the subset $R(x)$ is connected (in the sense  of Definition \ref{D75}).
\end{lem}

\begin{theorem}\label{th:D4KU1_KripCompl}
The logics ${\bf K4}+Ku$, ${\bf D4}+Ku$ are weakly canonical, and thus Kripke complete.
\end{theorem}

A proof of \ref{th:D4KU1_KripCompl} based on Lemma \ref{L77} and a 1-modal version of Lemma \ref{mcl} is straightforward, cf.
\cite{Sh90} or \cite{LB2} (the latter paper proves the same for $Ku_n$).

Hence we obtain
\begin{theorem}
The logic $\lgc{DT_1K}:=\lgc{DT_1}+Ku$ is weakly canonical, and thus  Kripke complete.
\end{theorem}

\begin{proof}

(Sketch.) For the axiom $Ku$ the argument from the proof of \ref{th:D4KU1_KripCompl} is still valid due to definability of all maximal clusters (Lemma \ref{mcl}). The remaining axioms are
Sahlqvist formulas.
{\hspace*{\fill} }\end{proof}

\begin{theorem}\label{th:DT_1CK_KripCompl}
The logic $\lgc{DT_1CK}:=\lgc{DT_1K}+AC^{\sharp u}$ is weakly canonical, and thus Kripke complete.
\end{theorem}

\begin{proof}
We can apply the previous theorem and
Lemma \ref{lem_KripkeCompletAC}.
{\hspace*{\fill} }\end{proof}

Completeness theorems from this section can be refined: in the next section we will prove the fmp for the logics considered above.

\section{The finite model property of ${\bf  D4K}$, $\lgc{DT_1K}$, and $\lgc{DT_1CK}$}

For the logic ${\bf  D4}+Ku$ the first proof of the fmp  was given in \cite{Sh90}. Another proof (also for ${\bf  D4}+Ku_n$) was proposed by M. Zakharyaschev \cite{Z}; it is based on a general and powerful method.

In this section we give a simplified version of the proof from \cite{Sh90}. It is based on a standard filtration method, and the same method is also applicable to 2-modal logics   $\lgc{DT_1K}$, $\lgc{DT_1CK}$.



\begin{theorem}\label{thm:fmp_DT_1CK}
The logics $\lgc{DT_1K}$ and $\lgc{DT_1CK}$ have the finite model property.
\end{theorem}

\begin{proof}
Let $\bLa$ be one of these logics.
Consider an $m$-formula $A\nin\Lambda$.
Take a generated submodel $M = (W, R, R_D, \vp)$ of the $m$-restricted canonical model of $\Lambda$ such that $M,u\nmo A$ for some $u$. As we know, its frame is basic and its $R$-maximal clusters are definable (Lemma \ref{mcl}).

Put
\begin{align*}
\Psi_0 &:= \setdef[\be(C)]{C\ \hbox{is an $R$-maximal cluster in $M$}}, \\
\Psi_1 &:= \set{A} \cup \setdef[\overline{\sqr}\gam]{\gam\ \hbox{is a Boolean combination of formulas from $\Psi_0$}}, \\
\Psi &:= \hbox{the closure of $\Psi_1$ under subformulas.}
\end{align*}

The set $\Psi$ is obviously finite up to equivalence in $\bLa$.

Take the filtration $M' = (W', R', R'_D, \vp')$ of $M$ through $\Psi$ as in Lemma \ref{LF1}. By that lemma, $F':=(W', R', R'_D)\mo {\bf KT_1}$.
The seriality of $R'$ easily follows from the seriality of $R$.

Next, if $\bLa=\lgc{DT_1CK}$, the frame $(W, R, R_D)$ is connected by Lemma \ref{lem_KripkeCompletAC}. So for any $x,y \in W$ there is an $R$-path from $x$ to $y$. $aRb$ implies $a^\sim R'b^\sim$, so
there is an $R'$-path from $x^\sim$ to $y^\sim$ in $F'$. Therefore
$F'\mo AC^{\sharp u}$.
It remains to show that $F'\mo Ku$.
Consider an $R'$-irreflexive point $x^\sim \in W'$ and assume that
$R' (x^\sim )$ is disconnected.
Let $V$ be a nonempty proper connected component of $R'(x^\sim)$. Consider
\begin{align*}
\Delta &:= \setdef[C]{\exists y (y^\sim \in V \;\&\; C\subseteq R(y) \;\&\;  C \hbox{ is an $R$-maximal cluster in $M$)}}; \\
B &:= \bigvee\limits_{C \in \Delta} \be (C),
\end{align*}
where $\be (C)$ is from Lemma \ref{mcl}.
Note that
\[
\leqno(1)\quad
z\in C\conj C\in\Del\Ra z^\sim \in V.
\]
In fact, if $C\in\Del$, then for some $y^\sim \in V$ we have $yRz$; hence
$y^\sim R'z^\sim$, so $z^\sim \in V$, by the  connectedness of $V$.

Let us show that for any $y^\sim \in R'(x^\sim)$
\[
\leqno(2)\quad
M', y^\sim \models B \ \ \hbox{iff}\ \  M, y \models B \ \ \hbox{iff}\ \  y^\sim \in V,
\]
i.e., $B$ defines $V$ in $ R'(x^\sim)$.

The first equivalence holds by the Filtration Lemma, since $B \in \Psi_1$.

Let us prove the second equivalence. To show `if', suppose $y^\sim\in V$.
By Lemma \ref{L82}, in the restricted canonical model there is a maximal cluster $C$ $R$-accessible from $y$; then $M, y \models \be(C)$. We have $C \in \Delta$, and thus $M, y \models B$.

To show `only if', suppose $y^\sim \not\in V$, but $M, y \models B$. Then $M,y \models \be (C)$, for some $C \in \Delta$,
hence $C \subseteq R(y)$, i.e., $yRz$ for some (and for all) $z\in C$; so it follows that $y^\sim R' z^\sim$.
Thus $y^\sim$ and $z^\sim$ are in the same connected component of $R'(x^\sim)$, which implies $z^\sim \nin V$. However,  $z^\sim \in V$ by (1), leading to a contradiction.

By Proposition \ref{cm} all substitution instances of $Ku$ are true in $M$. So
\[
M\mo Ku(B):=\sqr (\oB B \lor \oB\lnot B) \to \sqr B \lor \Box \lnot B.
\]

Consider an arbitrary $y\in R(x)$. Then for any $z\in R(y)$, $y^\sim$ and $z^\sim$ are in the same connected component of $R'(x^\sim)$.
Thus $y^\sim$ and $z^\sim$ are both either in $V$ or not in $V$, and so by (2), both of them satisfy either $B$ or $\lnot B$.
Hence $M, y\models \oB B \lor \oB\lnot B$. Therefore, $x$ satisfies the premise of $Ku(B)$. Consequently, $x$ must satisfy the conclusion of $Ku(B)$. Thus $M,x \models \sqr B$ or $M,x \models \sqr \lnot B$.
Since $\sqr B, \sqr \lnot B \in\Psi_1$, the Filtration Lemma implies
$M', x^\sim\mo \sqr B$ or $M', x^\sim\mo \sqr\neg B$. Eventually by (2),
$V=R'(x^\sim)$ or $V=\emp$, which
contradicts the assumption about $V$.

To conclude the proof, note that $A\in \Psi$, so by the Filtration Lemma $M',u^\sim\nvDash A$.
As we have proved, $F' \models \bLa$.
Therefore $\bLa$  has the fmp.
{\hspace*{\fill} }\end{proof}

\begin{theorem}\label{T101}
The logic ${\bf D4K}$  has the finite model property.
\end{theorem}

\begin{proof}
Use the argument from the proof of \ref{thm:fmp_DT_1CK} without the second relation.
{\hspace*{\fill} }\end{proof}

Thanks to the fmp, we have a convenient class of Kripke frames for the logic $\lgc{DT_1CK}$. This will allow us to prove the topological completeness result in the next section.

\section{The dd-logic of $\Real^n$, $n\ge 2$.}

This section contains the main result of the Chapter. The proof is based on the fmp theorem from the previous section and a technical construction of a dd-morphism presented in the Appendix.
	
In this section $\norm{\cdot}$ denotes the standard norm in $\Real^n$, i.e. for $x\in \Real^n$
\[
\norm{x} = \sqrt{x_1^2+\ldots+x_n^2}.
\]

We begin with some simple observations on connectedness. For a path $\alpha = w_0 w_1 \ldots
w_{n}$  in a $\lgc{K4}$-frame $(W,R)$ we use the notation $ \overline{R}(\alpha) :=
\bigcup\limits_{i=0}^n \overline{R} (w_{i})$. A path $\al$ is called {\em global} (in $F$) if
$\overline{R} (\al)=W$.

\begin{lem}\label{globpath}
Let $F = (W, R)$ be a finite connected $\bbK4$-frame, $w,v\in W$. Then there exists a global path from $w$ to $v$.
\end{lem}

\begin{proof}
In fact, in the finite connected
graph $(W, R^\pm)$ the vertices $w,v$ can be connected by a path visiting all the vertices (perhaps,
several times). {\hspace*{\fill} }
\end{proof}

\begin{lem}\label{LPF0}
Let $F = (W, R, R_D)$ be a finite rooted $\lgc{DT_1CK}$-frame. Then the set of all $R_D$-reflexive points in $F$ is connected.
\end{lem}

\begin{proof}
Let $x,y$ be two $R_D$-reflexive points. Since $(W,R)$ is connected, there exists a path connecting
$x$ and $y$. Consider such a path $\al$ with the minimal number $n$ of $R_D$-irreflexive points, and
let us show that $n=0$.

Suppose not. Take an $R_D$-irreflexive point $z$ in $\al$; then
$\al=x\ldots uzv\ldots y$, for some $u,v$, and it is clear that $zRu$, $zRv$, since $z$ is strictly $R$-minimal. By Lemma \ref{L77}, $R(z)$ is connected, so $u$, $v$ can be connected by a path $\be$ in $R(z)$. Thus in $\al$ we can replace the part $uzv$ with $\be$, and
the combined  path $x\ldots\be\ldots y$ contains $(n-1)$ $R_D$-irreflexive points, which contradicts the minimality of $n$.
{\hspace*{\fill} }
\end{proof}

\begin{lem}\label{lem:pathinaFrame}
Let $F = (W, R, R_D)$ be a finite rooted $\lgc{DT_1CK}$-frame and let $w', w'' \in W$ be
$R_D$-reflexive. Then there is a global path $\alpha = w_0 \ldots w_n$ in $(W,R)$ such that $w' =
w_0$, $w_n = w''$ and all $R_D$-irreflexive points occur only once in $\alpha$.
\end{lem}

\begin{proof}
Let $\set{u_1, \,\ldots, \;u_k}$ be the $R_D$-irreflexive points.
By connectedness there exists paths $\alpha_0$, \ldots, $\alpha_{k}$
respectively from $w'$ to $u_1$, from $u_1$ to $u_2$, \ldots, from $u_k$ to $w''$.

By Lemma \ref{LPF0}, the set $W' := W - \set{u_1, \ldots, u_k}$ is connected.
Hence we may assume that each $\alpha_i$ does not contain $R_D$-irreflexive points except its ends.
Also there exists a loop $\beta$ in $F':=F|W'$ from $w''$ to $w''$ such that
\hbox{$W - \bigcup\limits_{i=1}^{k-1} \overline{R}(\alpha_i) \subseteq \overline{R}(\beta).$}
\begin{figure}[h]
 \centering
 \includegraphics[width=0.85\textwidth]{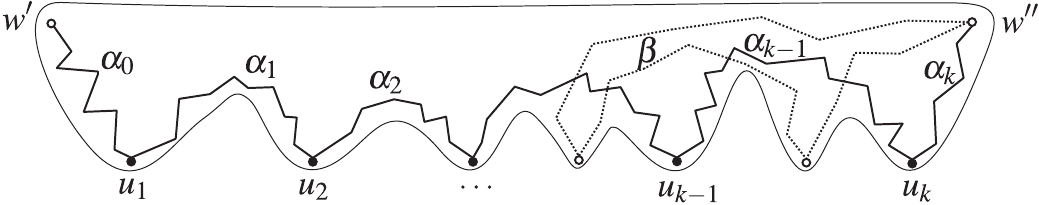}
 \caption{Path $\alpha$.}
 \label{fig:PathAlpha}
\end{figure}
Then we can define $\alpha$ as the joined path $\alpha_0\ldots \alpha_k
\beta$, (Fig. \ref{fig:PathAlpha}). {\hspace*{\fill} }\end{proof}

\begin{propo}\label{pr:dd-morphismRntoF}
For  a finite  rooted $\lgc{DT_1CK}$-frame $F = (W, R, R_D)$ and $R$-reflexive points $w', w'' \in
W$, the following holds.
\begin{description}
\item[(a)] If  $X=\set{ x\in \Rn  \mid  ||x|| \le r }$, $n \ge 2$,
then there exists $f: X \cpmor F$ such that $f(\partial X) = \set{w'}$;
\item[(b)] If $0\leq r_1<r_2$ and
\begin{align*}
 X& = \setdef[x\in \Rn]{r_1\leq ||x||\leq r_2 }, \\
Y'& =\setdef[x\in \Rn]{||x||=r_1 },\  Y'' =\setdef[x\in \Rn]{||x||=r_2 },
\end{align*}
then there exists $f: X \cpmor F$ such that
 $f(Y') = \set{w'}$, $f(Y'') = \set{w''}$.
\end{description}
\end{propo}

\begin{proof}
By induction on $\abs{W}$. Let us prove (a) first.
There are five cases:

\textbf{(a1)} $W = R(b)$ (and hence $b R b$) and $b = w'$.
Then there exists $f:X\dri^d (W,R)$.
In fact, let $C$ be the cluster of $b$ (as a subframe of $(W,R)$). Then $(W,R)=C$ or $(W,R)=C\cup F_1\cup\ldots\cup F_l$, where the
$F_i$ are generated by the successors of $C$. If $(W,R)=C$, we apply Proposition \ref{P54};
otherwise we apply Lemma \ref{refroot} and IH.

By \ref{pr:DS_AT1} it follows that $R_D$ is universal. And so by \ref{lem:dmortodDmor}(3) $f$ is a dd-morphism.

\medskip
\textbf{(a2)} $W = R(b)$ and not $w'Rb$. We may assume that $r=3$. Put
\[
 X_1 := \set{x  \mid ||x|| \le 1},~
 Y :=  \set{ x  \mid  1\le ||x|| \le  2 },~
 X_2 :=  \set{ x  \mid  2 \le ||x|| \le  3 }.
\]

By the case (a1), there is $f_1: X_1 \cpmor F$ with $f_1(\partial X_1) = \set{b}$. Let $C$ be a
maximal cluster in $R(w')$. By \ref{P54} there is $g: \BI Y \dri^d C$. Since $R(w')\neq W$,
we can apply IH to the frame $F':=F^{w'}_{\fo}$ and construct a dd-morphism $f_2:
X_2 \cpmor F'$ with $f_2(\partial X_2) = \set{w'}$. Now since $f_i(\partial X_i)\sbq R^{-1}(C)$, the
Glueing lemma \ref{Glue} is applicable. Thus $f:X\dri^d F$ for $f:=f_1\cup f_2\cup g$ (See Fig. \ref{fig:dd-morphism_a23}, Case (a2)). Note that
$\partial X\subset
\partial X_2$, so $f(\partial X)=f_2(\partial X)=\{w'\}$.

As in the case (a1), $f$
is a dd-morphism by \ref{lem:dmortodDmor}.

\medskip
\textbf{(a3)} $(W,R)$ is not rooted. By Lemma \ref{lem:pathinaFrame} there is a global path $\alpha$ in $F$ with a single occurrence of every $R_D$-irreflexive point. 
We may assume that $ \alpha =
b_0 c_0 b_1 c_1\ldots c_{m-1} b_m$, $b_m = w'$ and for any $i<m$ $c_i \in C_i \subseteq R(b_i) \cap R(b_{i+1})$, where $C_i$ is an $R$-maximal cluster.
Such a path is called \emph{reduced}.
For $0 \le j \le m$ we put $F_j := F|\overline{R}(b_j)$.

Since $(W,R)$ is not rooted, each $F_j$ is of smaller size than $F$, so we can apply the induction
hypothesis to $F_j$.
We may assume that
$$X = \set{ x \mid  ||x||\leq 2m+1}, \ Y=\set{ x \mid  ||x||=2m+1}.$$
Then put
\[
X_i := \set{ x \mid  ||x|| \leq  i+1  }\mbox{ for }0 \le i \le 2m,
\]
\[
Y_i := \partial X_i,~ \Del _i := \BC (X_{i} -X_{i-1})\mbox{ for }0 \le i \le 2m.
\]
By IH and Proposition \ref{P54} there exist
\begin{align*}
f_0:& X_0 \dri^{dd} F_0 \hbox{ such that } f_0(Y_0) =  \{c_0\},\\
f_{2j}:&  \Del _{2j} \dri^{dd}F_{j} \hbox{ such that } f_{2j}(Y_{2j}) = \{c_j\}, \ f_{2j}(Y_{2j-1}) = \{c_{j-1}\}\ \mbox{ for }1 \le j \le m,\\
f_{2j-1}:& \I \Del_{2j+1} \dpmor C_j \mbox{ for }0 \le j \le m-1.
\end{align*}

\begin{figure}
 \centering
\begin{tabular}{ccc}

 \includegraphics[width=0.3\textwidth]{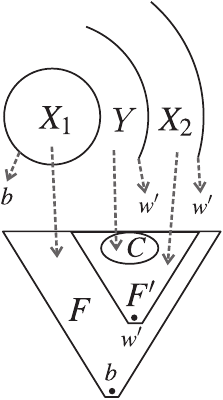}
&\qquad\qquad&
 \includegraphics[width=0.5\textwidth]{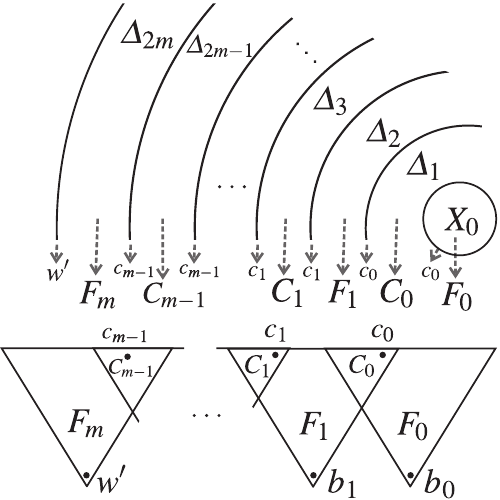}\\
Case (a2)&& Case (a3)
\end{tabular}
 \caption{dd-morphism $f$}
 \label{fig:dd-morphism_a23}
\end{figure}

One can check that $f:X\dri^{dd} F$ for $f:=\bigcup\limits_{j=0}^{2m}f_j$ (Fig. \ref{fig:dd-morphism_a23}).

\textbf{(a4)} $W = \overline{R}(b)$, $\lnot b R_D b$ (and so $\lnot b R b$).
We may assume that
$$X = \set{ x \mid  ||x||\leq 2}, \ Y=\set{ x \mid  ||x||=2}.$$
Then similar to case (a3) put
\[X_0 := X, ~Y_0 := Y, ~
X_i := \set{ x \mid  ||x|| \leq \frac 1i  },~ Y_i := \partial X_i,~ \Del _i := \BC (X_{i} -X_{i+1}), ~ (i>0).
\]
Consider the frame $F' := F|W'$, where $W' = W - \set{b}$. Note that $w'\in W'$, since $w'Rw'$, by
the assumption of \ref{pr:dd-morphismRntoF}. By Lemma \ref{L77} $F'$ is connected, and
thus $F'\mo \lgc{DT_1CK}$. By Lemma \ref{lem:pathinaFrame} there is a reduced global path $\alpha = a_1
\ldots a_m$ in $F'$ such that $ a_1 = w'$. Let
\[
\gamma = a_1 a_2 \ldots a_{m-1} a_m a_{m-1} \ldots a_2 a_1 a_2 \ldots
\]
be an infinite path shuttling back and forth through $\alpha$. Rename the  points in $\gamma$:
\begin{equation}
 \gamma =  b_0 c_0 b_1 c_1 \ldots b_m c_m b_{m+1} \ldots
\end{equation}
Again as in the case (a3) we put $F_j := F|\overline{R}(b_j)$, and assume that $c_j \in C_j$ and $C_j$ is an $R$-maximal cluster.
By IH there exist
\begin{align*}
f_0: & \Del _0 \dri^{dd} F_0 \hbox{ such that } f_0(Y_0) = \{b_0\} = \{w'\}, \ f_1(Y_1) =\{ c_0\},\\
f_{2j}: & \Del _{2j} \dri^{dd}F_{j} \hbox{ such that } f_{2j}(Y_{2j}) = \{c_{j-1}\}, \
f_{2j}(Y_{2j+1}) = \{c_{j}\}\hbox{ for }j>0,
\end{align*}
and by Proposition \ref{P54} there exist
$f_{2j+1}: \I \Del _{2j+1} \dpmor C_j$.
Put
\[
f(x) := \left\{
\begin{array}{ll}
b &\hbox{ if } x = {\bf 0},\\
f_{2j}(x)&\hbox{ if } x \in \Delta_{2j},\\
f_{2j+1}(x)&\hbox{ if } x \in \I\Delta_{2j+1},
\end{array}
\right.
\]
One can check that $f$ is d-morphic (Fig. \ref{fig:dd-morphism_a45}).
%
%
%
%
%
%
%

\medskip
\textbf{(a5)} $W = \overline{R}(b)$, $\lnot b R b$ and $b R_D b$. Then $R_D$ is universal,
$w'\ne b$. Put
\[
 X' :=  \set{x  \mid ||x|| < 1},  ~
X_4 := \set{ x \mid  1 \le ||x|| \le 2},
\]
and let $X_1,X_2$ be two disjoint closed balls in $X'$,
$X_3 := X'- X_1 - X_2$.

Let $C$ be a maximal cluster in $R(w')$,
$F' := F|R(w')$. Then there exist:
\begin{align*}
 f_i:& X_i \dri^d (W,R) \mbox{ for }i=1,2\mbox{ such that } \ f_i(\partial X_i) = \set{w'}, \hbox{ by the case (a4),}\\
 f_3:& X_3 \dri^d C, \hbox{ by Proposition \ref{P54},}\\
f_4:& X_4 \cpmor F' \mbox{ such that } \ f_4(\partial X_4) = \set{w'}, \hbox{ by the induction hypothesis.}
\end{align*}

\begin{figure}
 \centering
\begin{tabular}{ccc}

 \includegraphics[width=0.4\textwidth]{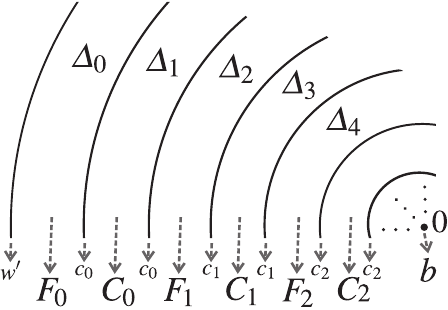}
&\qquad\qquad&
 \includegraphics[width=0.4\textwidth]{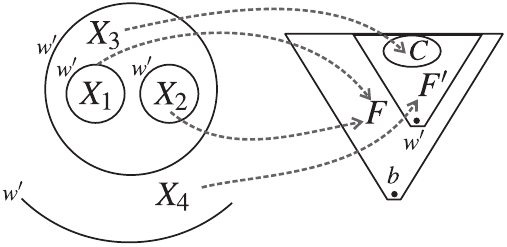}\\
Case (a4)&& Case (a5)
\end{tabular}
 \caption{dd-morphism $f$}
 \label{fig:dd-morphism_a45}
\end{figure}

Put $f:=f_1\cup f_2\cup f_3\cup f_4$ (Fig. \ref{fig:dd-morphism_a45}). Then $f(\partial \X) = \set{w'}$.

By Lemma \ref{Glue} (b) $f_1\cup f_2:X_1\cup X_2\dri^d F$, and hence $f:X\dri^d F$ by
Lemma \ref{Glue} (a).
 $f$ is manifold at $b$, thus it is a dd-morphism by
\ref{lem_pmorphism}.

\medskip
Now we prove (b). There are three cases.

\textbf{(b1)} $w' = w'' = b$ and $W = R(b)$.
The argument is the same as in the case (a1), using
Proposition \ref{P54},
Lemma \ref{refroot}, the induction hypothesis, and
Proposition \ref{lem:dmortodDmor}.

\medskip
\textbf{(b2)} $w' = w'' = b$, but $W \neq R(b)$. Consider a maximal cluster $C \subseteq R (b)$.
Since all spherical shells  for different $r_1$ and $r_2$ are homeomorphic, we assume that $r_1 = 1$, $r_2 = 4$. Consider the sets
\begin{equation*}
X_1 :=  \set{x \mid 1 \le ||x|| \le 2},\ \
 X' :=  \set{x \mid 2 < ||x|| < 3},\ \
X_3 :=  \set{x \mid 3 \le ||x|| \le 4},
\end{equation*}
and let $X_0\subset X'$ be a closed ball, $X_2 := X' - X_0$.
Let $F' :=F| R (b)$. There exist
\begin{align*}
 f_1:& X_1 \cpmor F'\mbox{ such that } f_1(\partial X_1) = \set{b}, \hbox{ by the case (b1)},\\
 f_2:& X_2 \dri^d C, \hbox{ by Proposition \ref{P54}},\\
 f_3:& X_3 \cpmor F'\mbox{ such that } \ f_3(\partial X_3) = \set{b}, \hbox{ by the case (b1)},\\
 f_0:& X_0 \cpmor F\mbox{ such that }  \ f_4(\partial X_0) = \set{b},
 \hbox{ by the statement (a) for }F.
\end{align*}
One can check that $f:X\cpmor F$ for $f:=f_0\cup f_1\cup f_2\cup f_3$.

%
%

\textbf{(b3)} $w' \ne w''$ and for some $b\in W$, $W = R(b)$, so $F$ has an $R$-reflexive root. Let
\[
F_1 := F|R(w'), ~ F_2 := F|R(w''),
\]
and let $C_i$ be an $R$-maximal cluster in $F_i$ for $i\in \set{1,2}$.

We assume that $r_1 = 1$, $r_2 = 6$ and consider the sets
\begin{align*}
X_i &:=  \set{x \mid i \le ||x|| \le i+1}, \  i \in \set{1, \ldots, 5}.
\end{align*}

By the case (b1) and Proposition \ref{P54} we have
\begin{align*}
 f_1:&~ X_1 \cpmor F_1 \mbox{ such that } f_1(\partial X_1) = \set{w'}, &f_2:&~ \I X_2 \dri^d C_1,\\
 f_3:&~ X_3 \cpmor F\mbox{ such that } \ f_3(\partial X_3) = \set{b},   &f_4:&~ \I X_4 \dri^d C_2, \\
 f_5:&~ X_5 \cpmor F_2 \mbox{ such that } f_1(\partial X_5) = \set{w''}.
\end{align*}

\begin{figure}
 \centering
\begin{tabular}{c|c}

 \includegraphics[width=0.45\textwidth]{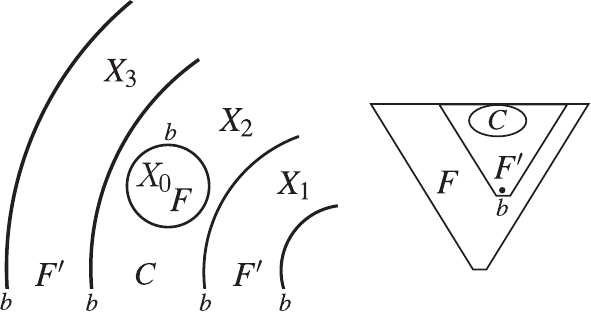}
\qquad&\qquad
 \includegraphics[width=0.3\textwidth]{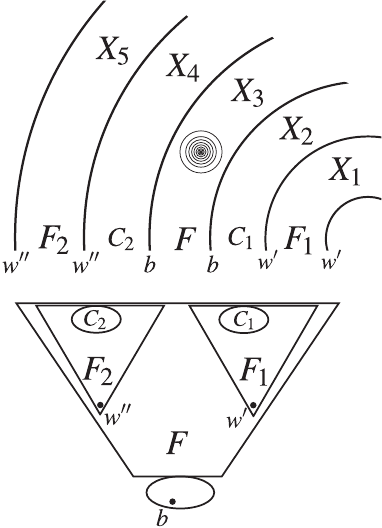}\\
Case (b2)& Case (b3)
\end{tabular}
 \caption{dd-morphism $f$}
 \label{fig:dd-morphism_b23}
\end{figure}

One can check that $f:X\cpmor F$ for $f:= \bigcup\limits_{i=1}^5 f_i$ (Fig. \ref{fig:dd-morphism_b23}, Case (b3)).



\medskip
\textbf{(b4)}  $w' \ne w''$ and $W \ne R(b)$ for any $b\in W$. By Lemma \ref{LPF0} there is a reduced path
$\alpha = b_0 c_0 b_1 \ldots c_{m-1} b_m$ from $b_0 = w'$ to $b_m = w''$ that does not contain $R_D$-irreflexive
points, $c_i \in C_i$, where $C_i$ is an $R$-maximal cluster. We may also assume that
\begin{equation}\label{eq:Rbi}
 \overline{R}(b_i) \ne W, \hbox{ for any $i\in \set{1, \ldots, m-1}$}.
\end{equation}
In fact, if the frame $(W,R)$ is not rooted, then (\ref{eq:Rbi}) obviously holds. If $(W,R)$
is rooted, then its root $r$ is irreflexive and by Lemma \ref{L77}, $R(r)$ is connected, so
there exists a path $\alpha$ in $R(r)$ satisfying (\ref{eq:Rbi}).
Put
\[
F_0 := F,~
F_j := F|R(b_j), 1 \le j \le m.
\]

Assuming that $r_1 = 1$, $r_2 = 2m+1$ we define
\begin{align*}
X_i &:= \setdef{||x|| \leq i +1 },~ Y_i := \partial X_i~ (\mbox{for }0 \le i \le 2m+1),\\
\Del _i &:= \BC (X_{i+1} -X_{i}) \ (\mbox{for }0\le i\le 2m).
\end{align*}

\begin{figure}[h]
 \centering
 \includegraphics[width=0.5\textwidth]{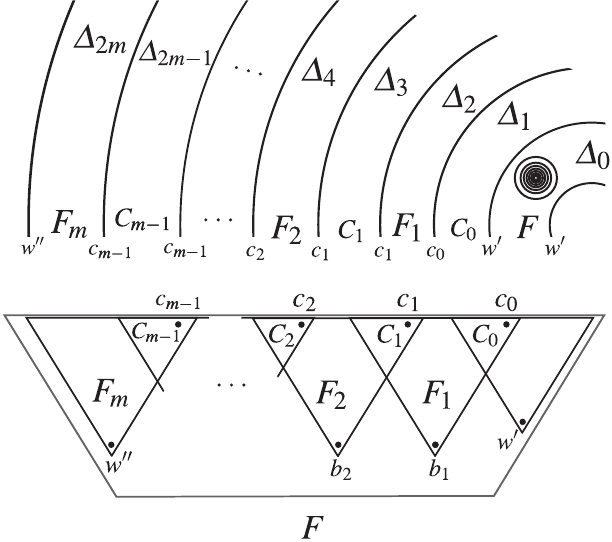}
 \caption{dd-morphism $f$, case (b4)}
 \label{fig:dd-morphism_b4}
\end{figure}

By the cases (b2), (b1), Proposition \ref{P54}, and the induction hypothesis there exist
\begin{align*}
&f_0:  \Del_0 \cpmor F=F_0 \hbox{ such that } f_0(Y_0) = f_0(Y_1)  = \set{w'}; \\
&f_{2j}:  \Del _{2j} \cpmor F_{j} \hbox{ such that } f_{2j}(Y_{2j+1}) = \set{c_{j}}, \ f_{2j}(Y_{2j}) = \{c_{j-1}\}\ (1 \le j \le m);\\
&f_{2j-1}:  \I \Del _{2j-1} \dri^d C_{j-1} \ (1 \le j \le m),\\
&f_{2m}:  \Del _{2m} \cpmor F_m \hbox{ such that } f_{2m}(Y_{2m}) = \set{c_{m}}, \ f_{2m}(Y_{2m+1}) =
\set{w''}.
\end{align*}

We claim that $f:X\cpmor F$ for $f:= \bigcup\limits_{i=0}^{2m} f_i$ (Fig. \ref{fig:dd-morphism_b4}).
First, we prove by induction using Lemma \ref{Glue} (see previous cases) that $f$ is a d-morphism.
Note that $f(Y')= f(Y_0) = \set{w'}$ and $f(Y'')= f(Y_{2m+1}) = \set{w''}$.

Second, there are no $R_D$-irreflexive points in $\al$, so all preimages of $R_D$-irreflexive points
are in $\Delta_0$; since $f_0$ is a dd-morphism, $f$ is 1-fold at any $R_D$-irreflexive point and
manifold at all the others. Thus $f$ is a dd-morphism by Proposition \ref{lem:dmortodDmor}.
\end{proof}

\begin{theorem}\label{T125}
 For $n \ge 2$, the dd-logic of $\BR^n$ is $\lgc{DT_1CK}$.
\end{theorem}

\begin{proof}
Since $\Rn$ is a locally 1-component connected dense-in-itself metric space, $\Rn \models^d
\lgc{DT_1CK}$.

Now consider a formula $A \notin \lgc{DT_1CK}$. Due to the fmp (Theorem \ref{thm:fmp_DT_1CK}) there
exists a finite rooted Kripke frame $F = (W, R, R_D) \mo \lgc{DT_1CK}$ such that $F \nvDash A$. By
Proposition $\ref{pr:dd-morphismRntoF}$ there exists $f:\Rn \cpmor F$. Hence $\Rn
\nvDash^d A$ by Lemma \ref{lem_pmorphism}. {\hspace*{\fill} }\end{proof}

\section{Concluding remarks}
\textbf{Hybrid logics.}
Logics with the difference modality are closely related to hybrid logics. The paper \cite{Litak06} describes a validity-preserving translation from the language with the topological and the difference modalities into the hybrid language with the topological modality, nominals and the universal modality.

Apparently a similar translation exists for dd-logics considered in our chapter. 
There may be an additional option --- to use `local nominals', propositional constants that may be true not in a single point, but in a discrete set. 
Perhaps one can also consider `one-dimensional nominals' naming `lines' or `curves' in the main topological space; there may be many other similar options.
\smallskip

\begin{figure}[h]
 \centering
 \includegraphics[width=4.1cm]{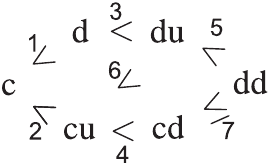}
 \caption{Correlation between topomodal languages.}
 \label{fig:languages}
\end{figure}

\textbf{Definability.}
Among several types of topological modal logics considered in this chapter dd-logics are the most
expressive. The correlation between all the types are shown in Fig. \ref{fig:languages}. 
A language $\cL_1$ is \emph{reducible} to $\cL_2$ ($\cL_1 \le \cL_2$) if every
$\cL_1$-definable class of spaces is
$\cL_2$-definable; $\cL_1<\cL_2$ if $\cL_1\le \cL_2$ and $\cL_2 \nleq \cL_1$. The non-strict reductions 1--7 in Fig. \ref{fig:languages}
are rather obvious. Let us explain, why 1--6 are strict.

The relations 1 and 2 are strict, since the c-logics of $\BR$ and $\BQ$ coincide \cite{MT}, while the cu- and d-logics are different \cite{Sh99, E1}.

The relation 3 is strict, since in d-logic without the universal modality we cannot express connectedness (this follows from \cite{E1}). The relations 4 and 6 are strict, since the cu-logics of $\BR$ and $\BR^2$ are the same \cite{Sh99}, while the cd- and du-logics are different \cite{Gab01, LB2}.

In cd- and dd-logic we can express \emph{global 1-componency}: the formula
\[
\DA (\overline{\sqr} p \lor \overline{\sqr} \lnot p) \to \DA p \lor \DA \lnot p
\]
is c-valid in a space $\X$ iff the complement of any point in $\X$ is connected. So we can distinguish the line $\Real$ and the circle ${\bf S^1}$. In du- (and cu-) logic
this is impossible, since there is a local homemorphism $f(t) = e^{it}$  from $\Real$ onto
$S^1$. It follows that the relation 5 is strict.
Our conjecture is that the relation 7 is strict as well.
\smallskip

\textbf{Axiomatization.}
There are several open questions about axiomatization and completeness of certain dd-logics.

1. The first group of questions is about the logic of $\BR$. On the one hand, in \cite{Kudinov08} it
was proved that $\mathbf{Lc}_{\ne}(\BR)$ is not finitely axiomatizable. Probably, the same method
can be applied to $\mathbf{Ld}_{\ne}(\BR)$. On the other hand,  $\mathbf{Lc}_{\ne}(\BR)$ has the fmp
\cite{Kud11}, and we hope that the same holds for the dd-logic. The decidability  of
$\mathbf{Ld}_{\ne}(\BR)$ follows from \cite{BG}, since this logic is a fragment of the universal
monadic theory of $\BR$; and by a result from \cite{Rey} it is PSPACE-complete. However, constructing
an explicit infinite axiomatization of $\mathbf{Lc}_{\ne}(\BR)$ or $\mathbf{Ld}_{\ne}(\BR)$ might be a serious technical problem.

2. A `natural' semantical characterization of the logic $\lgc{DT_1C}+Ku_2$ (which is a proper sublogic of $\mathbf{Ld}_{\ne}(\BR)$) is not quite clear. Our conjecture is that it is complete w.r.t. 2-dimensional cell complexes, or more exactly, adjunction spaces obtained from finite sets of 2-dimensional discs and 1-dimensional segments.

3. We do not know any syntactic description of dd-logics of 1-dimensional cell complexes (i.e., unions of finitely many segments in $\BR^3$ that may have only endpoints as common). 
 Their properties are probably similar to those  of $\mathbf{Ld}_{\ne}(\BR)$.

4. It may be interesting to study topological modal logics with the graded difference modalities $\DA_n A$ with the following semantics: $x\models \DA_n A$ iff there are at least $n$ points $y\ne x$ such that $y \models A$.

5. The papers \cite{MT} and \cite{Kremer} prove completeness and strong completeness of $\bbS4$ w.r.t.~any dense-in-itself metric space. The corresponding result for d-logics is completeness of $\mathbf{D4}$ w.r.t. an arbitrary dense-in-itself separable metric space. Is separabilty essential here? Does strong completeness hold in this case? Similar questions make sense for dd-logics. 

6. \cite{Gab01} presents a 2-modal formula cd-valid exactly in $T_0$-spaces. However, the cd-logic (and the dd-logic) of the class of $T_0$-spaces is still unknown. Note that the d-logic of this class has been axiomatized in \cite{BEG2011}; probably the same technique is applicable to cd- and dd-logics.

7. In footnote 7 we have mentioned that there is a gap in the paper \cite{Sh99}. Still we can prove that for any connected, locally connected metric space $\X$ such that the boundary of any ball is nowhere dense, $\BL\Bc_\fo(\X)=\mathbf{S4U} + AC$. But for an arbitrary connected metric space $\X$ we do not even know if $\BL\Bc_\fo(\X)$ is finitely axiomatizable.

8. Is it possible to characterize finitely axiomatizable dd-logics that are complete w.r.t.~Hausdorff spaces? metric spaces? Does there exist a dd-logic complete w.r.t. Hausdorff spaces, but incomplete w.r.t. metric spaces?

9. Suppose we have a c-complete modal logic $L$, and let ${\cal K}$ be the class of all topological spaces where $L$ is valid.  Is it always true that ${\bf Lc}_\fo({\cal K})=LU$? and ${\bf Lc}_{\ne}({\cal K})=LD$? 
Similar questions can be formulated for d-complete modal logics and their du- and dd-extensions.

10. An interesting topic not addressed in this chapter is the complexity of topomodal logics. In particular, the complexity is unknown for 
the d-logic (and the dd-logic) of $\BR^n$ ($n>1$).

We would like to thank anonymous referee who helped us improve the first version of
the manuscript.

The work on this chapter was supported by RFBR grants 11-01-00281-a,
11-01-00958-a, 11-01-93107-CNRS-a and the Russian President's grant 
NSh-5593.2012.1.


\end{document}